\documentclass[11pt]{article}
\usepackage{amssymb,amsmath}

\usepackage[dvips]{epsfig}
\usepackage{graphicx}

\newtheorem{Lemma}{Lemma}[section]
\newtheorem{Theorem}{Theorem}
\newtheorem{Proposition}[Lemma]{Proposition}
\newtheorem{Corollary}[Lemma]{Corollary}
\newtheorem{Remark}[Lemma]{Remark}
\newtheorem{Remarks}[Lemma]{Remarks}

\newtheorem{Hypothesis}[Lemma]{Hypothesis}

\newenvironment{Proof}
 {\begin{trivlist} \item[]{\bf Proof. }}
 {\hspace*{\fill}$\rule{.3\baselineskip}{.35\baselineskip}$\end{trivlist}}

\makeatletter\@addtoreset{figure}{section}\makeatother

\makeatletter\@addtoreset{equation}{section}\makeatother

\newcommand{\C}{\mathbb{C}}
\newcommand{\N}{\mathbb{N}}
\newcommand{\R}{\mathbb{R}}
\newcommand{\Z}{\mathbb{Z}}

\def\Re{\mathop{\mathrm{Re}}}

\def\Tr{\mathrm{Tr}}
\def\Det{\mathrm{Det}}

\newcommand{\id}{\mathrm{\,id}\,}
\newcommand{\rmO}{\mathrm{O}}
\newcommand{\rmo}{\mathrm{o}}
\newcommand{\rmd}{\mathrm{d}}
\newcommand{\rme}{\mathrm{e}}
\newcommand{\rmi}{\mathrm{i}}

\begin{document}

\begin{center}

\vspace*{4mm}

{\LARGE Diffusive stability of oscillations in reaction-diffusion systems}
\\[8mm]
\begin{minipage}[t]{0.4\textwidth}
\begin{center}
{\large Thierry Gallay}\\
Universit\'e de Grenoble I\\
Institut Fourier, UMR CNRS 5582\\
BP 74\\
38402 Saint-Martin-d'H\`eres, France
\end{center}
\end{minipage}
\begin{minipage}[t]{0.4\textwidth}
\begin{center}
{\large Arnd Scheel}
\\
University of Minnesota\\
School of Mathematics\\
206 Church St. S.E.\\
Minneapolis, MN 55455, USA
\end{center}
\end{minipage}
\end{center}

\vspace*{1cm}

\begin{abstract}
We study nonlinear stability of spatially homogeneous oscillations
in reaction-diffusion systems. Assuming absence of unstable linear
modes and linear diffusive behavior for the neutral phase, we prove
that spatially localized perturbations decay algebraically with the
diffusive rate $t^{-n/2}$ in space dimension $n$. We also compute
the leading order term in the asymptotic expansion of the solution,
and show that it corresponds to a spatially localized modulation of
the phase. Our approach is based on a normal form transformation in
the kinetics ODE which partially decouples the phase equation, at
the expense of making the whole system quasilinear. Stability is
then obtained by a global fixed point argument in temporally
weighted Sobolev spaces.
\end{abstract}

\thispagestyle{empty}

\vfill

{\bf Corresponding author:} Arnd Scheel


\newpage

\section{Introduction and main results}
\label{s:1}
Synchronization of spatially distributed dissipative oscillators has
been observed in a wide variety of physical systems. We mention
synchronization in yeast cell populations \cite{GCP}, fireflies
\cite{BB}, coupled laser arrays \cite{RTF}, and spatially homogeneous
oscillations in reaction-diffusion systems such as the
Belousov-Zhabotinsky reaction \cite{Win} and the NO+CO-reaction on a
Pt(100) surface \cite{VMMI}. Synchronization strikes us most when the
system size is large or the coupling strength is weak. Both situations
relate in natural ways to the regime of large Reynolds number in fluid
experiments, where one expects turbulent, incoherent rather than
laminar, synchronized behavior. Still, one finds synchronization as a
quite common, universal phenomenon, even in very large systems.

The aim of this article is to elucidate the robustness of spatially
homogeneous temporal oscillations in spatially extended systems,
under most general assumptions, without detailed knowledge of internal
oscillator dynamics or coupling mechanisms. In fact, quantitative
models are very rarely available for the systems mentioned above.
Instead, we make phenomenological assumptions, related to the
existence of oscillations and the absence of strongly unstable modes.
These assumptions typically guarantee asymptotic stability of a
spatially homogeneous oscillation in any \emph{finite-size} system,
when equipped with compatible (say, Neumann) boundary conditions. The
results in this article are concerned with \emph{infinite-size},
reaction-diffusion systems,
\begin{equation}\label{rdsystem}
  u_t \,=\, D\Delta u + f(u),\qquad u = u(t,x) \in\R^N\,, 
  \quad x\in\R^n\,, \quad t \ge 0\,,
\end{equation}
with positive coupling matrix $D \in \mathcal{M}_{N\times N}(\R)$, $D=D^T>0$,
and smooth kinetics $f\in C^\infty(\R^N,\R^N)$. In this spatially
continuous setup, working in the whole space $\R^n$ is an idealization 
which corresponds to the limit of small coupling matrix and/or
large domain size. We will briefly comment on the relation
between our results in the whole space and the stability of 
temporal oscillations in finite domains, below.  

To be specific, we make the following assumptions on the kinetics $f$
and the coupling matrix $D$.

\begin{Hypothesis}[Oscillatory kinetics]\label{h:1}
We suppose that the ODE $u_t = f(u)$ possesses a periodic solution 
$u_*(t) = u_*(t+T)$ with minimal period $T>0$. 
\end{Hypothesis}

In particular, to avoid trivial situations, we assume that the 
periodic orbit is not reduced to a single equilibrium. As is 
well-known, this is possible only if $N \ge 2$, i.e. if the 
system~\eqref{rdsystem} does not reduce to a scalar equation. 

In addition to existence we will make a number of assumptions on
the Floquet exponents of the linearized equation
\begin{equation}\label{e:lin}
  u_t \,=\, D\Delta u + f'(u_*(t))u\,,
\end{equation}
which is formally equivalent to the family of ordinary differential 
equations
\begin{equation}\label{e:linf}
  u_t \,=\, -k^2D u + f'(u_*(t))u\,, \quad k\in\R^n\,.
\end{equation}
For each fixed $k$ we denote by $F_k(t,s)$ the two-parameter
evolution operator associated to the linear time-periodic 
system \eqref{e:linf}, so that $u(t) = F_k(t,s)u(s)$ for any 
$t \ge s$. The asymptotic behavior of the solutions of \eqref{e:linf}
is well characterized by the Floquet multipliers of the system, 
that is, the eigenvalues of the period map $F_k(T,0)$. We 
shall rather work with the Floquet exponents $\lambda_1(k),\dots,
\lambda_N(k) \in \C/\rmi\omega\Z$, where $\omega = 2\pi/T$,   
which satisfy 
\[
  \mathrm{det}\,\Bigl(F_k(T,0) - \rme^{\lambda_j(k)T}\Bigr) 
  \,=\,0\,, \quad j = 1,\dots,N\,.
\]
Remark that any Floquet exponent which is simple (that is, of 
algebraic multiplicity one) depends smoothly on the parameter $k$. 
Also note that $\lambda_1 = 0$ is always a Floquet exponent for $k=0$. 
We refer to the set of Floquet exponents as the Floquet spectrum.

\begin{Hypothesis}[Marginally stable spectrum]
\label{h:2}
We suppose that the Floquet spectrum in the closed half-space 
$\{\Re\lambda\ge 0\}$ is minimal. More precisely, we assume that\\
(i) The Floquet spectrum in the closed half-space $\{\Re\lambda\ge
0\}$ is nonempty only for $k = 0$, in which\\ 
\null\hskip 18pt case it consists of a simple Floquet exponent 
$\lambda_1 = 0$;\\
(ii) Near $k = 0$, the neutral Floquet exponent continues as 
$\lambda_1(k) = -d_0 k^2+\rmO(k^4)$ for some $d_0>0$.
\end{Hypothesis}

We emphasize that these assumptions are satisfied for an open class of
reaction-diffusion systems. In particular the expansion (ii) with some
$d_0\in\R$ is a consequence of the simplicity of the Floquet exponent
$\lambda_1 = 0$ at $k=0$, and of the symmetry $k\mapsto -k$; assuming
$d_0>0$ is therefore robust. In fact, it is not difficult to show that
\begin{equation}\label{e:d0def}
  d_0 \,=\, \frac{\int_0^T (U_*(t),Du_*'(t))\,\rmd t}{\int_0^T
    (U_*(t),u_*'(t))\,\rmd t}\,,
\end{equation}
where $(\cdot,\cdot)$ denotes the usual scalar product in $\R^N$, 
and $U_*(t)$ is the (unique nontrivial) bounded solution of the 
adjoint equation
\begin{equation}\label{e:adj}
  -U_t \,=\, f'(u_*(t))^T U\,. 
\end{equation}
Of course, a necessary condition for Hypothesis~\ref{h:2} to hold 
is that $u_*(t)$ be a {\em stable} periodic solution of the ODE 
$u_t = f(u)$, but this assumption alone is not sufficient in general, 
except if the diffusion matrix is a multiple of the identity.
Indeed, even if $N = 2$, one can find examples of periodic solutions
which are asymptotically stable for the ODE dynamics, but become 
unstable if a suitable diffusion is added \cite{Ri1,Ri2,Ri3}. 
One possible scenario, which is usually called {\em phase instability}
or {\em sideband instability}, is that the coefficient $d_0$ 
be negative, in which case the periodic orbit is unstable with 
respect to long-wavelength perturbations. It may also happen that
the Floquet spectrum is stable for $k$ in a neighborhood of the 
origin, but that there exists an unstable Floquet exponent for
some $k_* \neq 0$, and therefore for all $k$ in a neighborhood of
$k_*$. This mechanism is reminiscent of the {\em Turing instability}
for spatially homogeneous equilibria. In Section~\ref{s:ex} below, 
we give an example of a simple 2-species system which exhibits 
both kinds of instabilities depending on the choice of the 
parameters. 

In order to state our results, we introduce a function space which
measures both the spatial localization and the amplitude of the
perturbation to our spatially homogeneous profile $u_*$. We will
consider initial perturbations in the space of functions $X =
L^1(\R^n) \cap L^\infty(\R^n)$, with target space $\R^N$, equipped
with the norm
\[
  \|v\|_X \,=\,  \int_{\R^n} |v(x)|\,\rmd x \,+\, \sup_{x \in 
  \R^n} |v(x)|\,,
\]
where ``sup'' here refers to the essential supremum. We will measure 
decay in the space $L^\infty(\R^n)$. Our first result is:

\begin{Theorem}\label{t:1}
Consider a reaction-diffusion system \eqref{rdsystem} on $\R^n$
with oscillatory kinetics, Hypothesis~\ref{h:1}, and marginally 
stable spectrum, Hypothesis~\ref{h:2}. Then there are $C,\delta>0$
such that for any initial data $u(0,x) = u_*(t_0) + v_0(x)$ with
$t_0\in\R$ arbitrary and $\|v_0\|_X \le \delta$, there exists a 
unique, smooth global solution $u(t,x)$ of \eqref{rdsystem} 
for $t \ge 0$. Moreover $u(t,x)$ converges to the periodic solution 
$u_*$ in the sense that
\begin{equation}\label{e:t1decay}
  \sup_{x \in \R^n}\Big|u(t,x)-u_*(t_0+t)\Big| \,\le\, \frac{C 
  \|v_0\|_X}{(1+t)^{n/2}}\,, \quad \hbox{for all }t \ge 0\,.
\end{equation}
\end{Theorem}

We emphasize that the perturbations we consider are localized in 
space, and therefore do not alter the overall phase $t_0$ of the 
periodic solution. We refer however to Section~\ref{s:ex} for a 
discussion of possible stability results in more general situations.
It is not difficult to verify that the decay rate in \eqref{e:t1decay}
is optimal. In fact, under the assumptions of Theorem~\ref{t:1}, 
one can even compute the leading term in the asymptotic expansion of the 
perturbation as $t \to +\infty$. Let 
\begin{equation}\label{e:GGdef}
  G(x) \,=\, \frac{1}{(4\pi d_0)^{n/2}}\,\exp\Bigl(
  -\frac{|x|^2}{4d_0}\Bigr)\,, \quad x \in \R^n\,,
\end{equation}
where $d_0 > 0$ is defined in \eqref{e:d0def}. Our second result is:

\begin{Theorem}\label{t:2}
Under the assumptions of Theorem~\ref{t:1}, the solution $u(t,x)$
of \eqref{rdsystem} can be decomposed as
\begin{equation}\label{e:udecomp}
  u(t,x) \,=\, u_*(t_0+t) + u_*'(t_0+t) \alpha(t,x) + \beta(t,x)\,, 
  \quad x \in \R^n\,, \quad t \ge 0\,, 
\end{equation}
where $\alpha : \R_+ \times \R^n \to \R$ and $\beta : \R_+ \times 
\R^n \to \R^N$ satisfy
\begin{eqnarray}\label{e:betaest}
  \|\beta(t,\cdot)\|_{L^1} + (1+t)^{n/2}\|\beta(t,\cdot)\|_{L^\infty}
  & \xrightarrow[t\to +\infty]{}& 0\,, \\ \label{e:alphaest}
  \|t^{n/2}\alpha(t,x\sqrt{t}) - \alpha_* G\|_{L^1 \cap 
  L^\infty} & \xrightarrow[t\to +\infty]{}& 0\,,
\end{eqnarray}
for some $\alpha_* \in \R$. In addition, 
\begin{equation}\label{e:alphastardef}
  \alpha_* \,=\, \int_{\R^n} \frac{(U_*(t_0),v_0(x))}{(U_*(t_0),u_*'(t_0))}
  \,\rmd x + \rmO(\delta^2)\,,
\end{equation}
where $U_*(t)$ is the bounded solution of the adjoint equation
\eqref{e:adj}. 
\end{Theorem}

In other words, the solution $u(t,x)$ of \eqref{rdsystem} satisfies
\begin{eqnarray*}
  u(t,x) &=& u_*(t_0+t) + u_*'(t_0+t)\,\frac{\alpha_*}{(4\pi d_0 t)^{n/2}}
  \,\rme^{-|x|^2/(4d_0t)} + \rmo(t^{-n/2}) \\
  &=& u_*\Bigl(t_0+t + \frac{\alpha_*}{(4\pi d_0 t)^{n/2}}
  \,\rme^{-|x|^2/(4d_0t)}\Bigr) + \rmo(t^{-n/2})\,, \quad \hbox{as }
  t \to +\infty\,.
\end{eqnarray*}
To leading order, the effect of the perturbation is thus a spatially
localized modulation of the phase of the periodic solution. As is
clear from the proof, the left-hand side of \eqref{e:betaest} decays
at least like $t^{-\gamma}$ as $t \to \infty$, for some $\gamma > 0$.
However, to specify a convergence rate in \eqref{e:alphaest}, it is
necessary to restrict ourselves to more localized perturbations. For
instance, if we assume in addition that $(1+|x|)v_0 \in L^1(\R^n)$,
then we can prove that the left-hand side of \eqref{e:alphaest} is
$\rmO(t^{-1/2})$ as $t \to \infty$.

Under Hypothesis~\ref{h:2}, if we consider system \eqref{rdsystem} in
a large bounded domain (say $x \in \Omega/\varepsilon$ where $\Omega
\subset \R^n$ is bounded) with Neumann boundary conditions, the
perturbations of the periodic solution $u_*(t)$ decay exponentially
\cite{henry}: if $\|v_0\|_{L^\infty} \le \delta$ for some small
$\delta > 0$, then $\|v(t)\|_{L^\infty} \le C \|v_0\|_{L^\infty}
\,\rme^{-\mu t}$ for some $\mu > 0$. However, the relaxation rate
$\mu$ and the size of admissible perturbations $\delta$ both depend on
$\varepsilon$, with typical scalings $\mu,\delta =
\rmO(\varepsilon^2)$ predicted by the spectral gap of the Laplacian on
the domain $\Omega/\varepsilon$.  This gap vanishes in the limit
$\varepsilon \to 0$, and Theorem~\ref{t:1} shows that exponential
decay is replaced by diffusive decay. Nevertheless, we expect our
results to give an accurate description of the intermediate
asymptotics for large bounded domains, if the initial perturbations
are sufficiently localized.

The type of diffusive decay that we establish in Theorems~\ref{t:1}
and \ref{t:2} has been observed in many other contexts. For instance,
localized perturbations of spatially periodic, stationary patterns in
the Ginzburg-Landau or Swift-Hohenberg equation exhibit a similar
diffusive behavior \cite{BK92,CE92,CEE92,Ka94,Sc96,Ue99}. At a
technical level, the approach in \cite{BK92,Sc96,Ue99} is based on
{\em renormalization group} theory, see for instance \cite{BK94}.
Roughly speaking, the method relies on the fact that the time-$T$ map
for the evolution of the perturbations becomes a contraction in a
space of localized functions when composed with an appropriate
rescaling, except for a neutral direction which specifies the profile
of the self-similar solution describing the leading order asymptotics.
A nice feature of renormalization theory is that it allows to
determine easily which terms in the nonlinearity are ``relevant''
(that is, potentially dangerous) for the stability analysis.  For
example, if we consider the nonlinear heat equation $u_t = \Delta u +
|u|^p$ in $\R^n$, with small and localized initial data, it is
well-known that the nonlinearity $|u|^p$ will not influence the decay
predicted by the linear evolution if $p > 1+2/n$, whereas
instabilities and even blow-up phenomena can occur if $p \le 1+2/n$
\cite{Fu66,CEE92}. In particular, quadratic terms (which arise
naturally in the Taylor expansion of any smooth function) are
``irrelevant'' if $n \ge 3$ and ``relevant'' if $n = 1,2$. For this
reason, diffusive stability is often easier to establish in high space
dimensions, when diffusion is strong enough to control all possible
nonlinear terms, whereas serious problems can occur in low dimensions.
This is the case in particular in the stability analysis of
one-dimensional spatially periodic patterns \cite{Sc96,EWW97}, where a
key step of the proof is to show that relevant ``self-coupling'' terms
actually do not occur in the evolution equation for the neutral
translational mode.

As one may expect from the discussion above, Theorem~\ref{t:1} is
rather easy to prove when $n \ge 3$. For completeness, we first settle
this case in Section~\ref{s:high} and then focus on the more
interesting situation where $n = 1$ or $2$. Here the idea is to
construct a {\em normal form} transformation for the ODE dynamics
which removes all ``relevant'' terms in the nonlinear PDE satisfied by
the perturbation. In Section~\ref{s:nf}, we show that this is
possible, at the expense of transforming the semilinear equation
\eqref{rdsystem} into a quasilinear parabolic system.  The next
important step is to obtain optimal decay estimates for the solutions
of the linearized perturbation equation, including maximal regularity
estimates, using the spectral assumptions in Hypothesis~\ref{h:2}.
Since the perturbation equation is translation invariant in space and
periodic in time, such bounds are relatively straightforward to obtain
via Fourier analysis, see Section~\ref{s:le}. Using these linear
estimates, we give in Section~\ref{s:nl} a proof of Theorem~\ref{t:1}
which is valid for $n \le 3$, hence covering the missing cases $n =
1,2$.  Instead of renormalization group theory, we prefer using a
global fixed point argument in temporally weighted spaces, as in
\cite{CE92,CEE92,Ka94}. After stability has been established, a rather
classical procedure, which is recalled in Section~\ref{s:asym}, allows
to derive the first-order asymptotics and to prove Theorem~\ref{t:2}
at least for $n \le 3$ (the higher dimensional case is again easier,
and left to the reader). In the final Section~\ref{s:ex}, we
illustrate our spectral assumptions on a simple, explicit example, and
we conclude with a short discussion including possible extensions of
our results.

\noindent
\textbf{Acknowledgments.} The authors thank Sylvie Monniaux
for useful discussions on maximal regularity in parabolic equations.
A.S. would like to thank the Universit\'e de Franche-Comt\'e for
generous support and hospitality during his stay, where part of this
project was carried out. A.S also acknowledges partial support by the
NSF through grant DMS-0504271.

\section{Stability in high dimensions}\label{s:high}

In this section we explore a straightforward and somewhat naive
approach to the stability of the periodic orbit $u_*(t)$ as a solution
of the reaction-diffusion system \eqref{rdsystem}.  This method gives
a simple proof of Theorem~\ref{t:1} in the high-dimensional case $n
\ge 3$, the main ingredient of which is an $L^p$-$L^q$ estimate for
the linearized evolution operator which will be established in
Section~\ref{s:le}. Without loss of generality, we assume from 
now on that the parameter $t_0$ in Theorems~\ref{t:1} and \ref{t:2}
is equal to zero (this is just an appropriate choice of the origin 
of time). 

Consider a solution $u(t,x) = u_*(t) + v(t,x)$ of \eqref{rdsystem}.
The perturbation $v$ satisfies the equation
\begin{equation}\label{e:pertv}
  v_t \,=\, D\Delta v + f'(u_*(t))v + N(u_*(t),v)\,, 
\end{equation}
where
\[
  N(u_*(t),v) \,=\, f(u_*(t)+v)-f(u_*(t))-f'(u_*(t))v\,. 
\]
The Cauchy problem for the semi-linear parabolic system
\eqref{e:pertv} is locally well-posed in the space $X = L^1(\R^n) 
\cap L^\infty(\R^n)$, see e.g. \cite{henry,pazy}. More precisely, 
for any $v_0 \in X$, there exists a time $\tilde T > 0$ (depending
only on $\|v_0\|_X$) such that \eqref{e:pertv} has a unique (mild) 
solution $v \in C^0([0,\tilde T],L^1(\R^n)) \cap C^0_b((0,\tilde T],
L^\infty(\R^n))$ satisfying $v(0) = v_0$. 

\begin{Remark}\label{r:H2data}
Due to parabolic regularization, the solution $v(t,x)$ of 
\eqref{e:pertv} is smooth for $t > 0$. For instance, there
exists $C > 0$ such that $\|v(t)\|_{H^2} \le C t^{-1}\|v_0\|_X$
for all $t \in (0,\tilde T]$. Therefore, in the proof of 
Theorem~\ref{t:1}, we can assume without loss of generality 
that the initial perturbation is small in the space $X \cap 
H^2(\R^n)$. 
\end{Remark}

To investigate the long-time behavior of the solutions of
\eqref{e:pertv}, we consider the corresponding integral equation
\begin{equation}\label{e:pertvint}
  v(t) \,=\, \mathcal{F}(t,0)v_0 + \int_0^t \mathcal{F}(t,s)
  N(u_*(s),v(s))\,\rmd s\,,
\end{equation}
where $\mathcal{F}(t,s)$ is the two-parameter semigroup associated
to the linearized equation \eqref{e:lin}. Due to our spectral 
assumptions (Hypothesis~\ref{h:2}), the operator $\mathcal{F}(t,s)$
satisfies the same $L^p$--$L^q$ estimates as the heat semigroup 
$\rme^{(t-s)\Delta}$. More precisely, anticipating the results of 
Section~\ref{s:le}, we have:

\begin{Proposition}[$L^p$--$L^q$ estimates]\label{p:10}
There exists a positive constant $C$ such that, for all $t > s$ 
and all $1 \le p \le q \le \infty$, we have
\begin{equation}\label{e:newlplq}
  \|\mathcal{F}(t,s)v\|_{L^q(\R^n)} \,\le\, \frac{C}{(t-s)^{\frac{n}2
  (\frac1p-\frac1q)}}\ \|v\|_{L^p(\R^n)}\,.
\end{equation}
\end{Proposition}

\begin{Proof}
The proof follows exactly the same lines as in Propositions~\ref{p:1},
\ref{p:2} and \ref{p:3}.  
\end{Proof}

By construction, the nonlinearity $N(u_*,v)$ in \eqref{e:pertv} 
is at least quadratic in $v$ in a neighborhood of the origin. 
More precisely, there exists a nondecreasing function $K : \R_+ 
\to \R_+$ such that, for all $t \in [0,T]$, 
\[
  |N(u_*(t),v)| \,\le\, K(R)|v|^2 \quad \hbox{whenever }
  |v| \le R\,.
\]
As was mentioned in the introduction, if the space dimension $n$ is 
greater or equal to $3$, the diffusive effect described in 
\eqref{e:newlplq} is strong enough to kill the potential instabilities 
due to the nonlinearity. In that case, nonlinear stability can 
therefore be established by a classical argument, which we briefly
reproduce here for the reader's convenience.

\medskip\noindent
\textbf{Proof of Theorem~\ref{t:1}} ($n \ge 3$).  
Fix $v_0 \in X = L^1(\R^n) \cap L^\infty(\R^n)$, and let $v \in
C^0([0,T_*),L^1(\R^n)) \cap C^0((0,T_*),L^\infty(\R^n))$ be the
maximal solution of \eqref{e:pertv} with initial data $v_0$. 
For $t \in [0,T_*)$ we denote
\[
  \phi(t) \,=\, \sup_{0 \le s \le t}\|v(s)\|_{L^1} + 
  \sup_{0 \le s \le t}(1+s)^{n/2}\|v(s)\|_{L^\infty}\,.  
\]
Using the integral equation \eqref{e:pertvint} and the linear
estimates \eqref{e:newlplq}, we easily find
\begin{align*}
\|v(t)\|_{L^1} &\le \|\mathcal{F}(t,0)v_0\|_{L^1} + 
  \int_0^t \|\mathcal{F}(t,s) N(u_*(s),v(s))\|_{L^1}\,\rmd s\\
&\le C\|v_0\|_{L^1} + CK(\phi(t))\int_0^t \|v(s)\|_{L^1}
  \|v(s)\|_{L^\infty} \,\rmd s\\
&\le C\|v_0\|_{L^1} + CK(\phi(t))\phi(t)^2 \int_0^t
 \frac{1}{(1+s)^{n/2}}\,\rmd s\,.
\end{align*}
Similarly, if $0 < t < 1$, we have
\begin{align*}
(1+t)^{n/2}\|v(t)\|_{L^\infty} &\le C\|v_0\|_{L^\infty} 
  + C\int_0^t \|N(u_*(s),v(s))\|_{L^\infty}\,\rmd s\\
&\le C\|v_0\|_{L^\infty} + CK(\phi(t))\phi(t)^2 \int_0^t 
  \frac{1}{(1+s)^n}\,\rmd s\,,
\end{align*}
while for $t \ge 1$ we can bound
\begin{align*}
(1+t)^{n/2}\|v(t)\|_{L^\infty} &\le C\|v_0\|_{L^1} + 
  C(1+t)^{n/2}\int_0^{t/2}\frac{1}{(t-s)^{n/2}}\ \|N(u_*(s),
  v(s))\|_{L^1}\,\rmd s\\
&\hspace{2cm} + C(1+t)^{n/2}\int_{t/2}^t \|N(u_*(s),v(s))\|_{L^\infty}
\,\rmd s\\
&\le C\|v_0\|_{L^1} + CK(\phi(t))\phi(t)^2\int_0^t\frac{1}{(1+s)^{n/2}}
\,\rmd s\,.  
\end{align*}
Now, since $n \ge 3$, we have $\int_0^\infty (1+s)^{-n/2}\,\rmd s
< \infty$ and we see that there exist positive constants $C_1, C_2$ 
(independent of $T_*$) such that
\begin{equation}\label{e:phibd}
  \phi(t) \,\le\, C_1 \|v_0\|_X + C_2 K(\phi(t))\phi(t)^2, \quad 
  \hbox{for all } t \in [0,T_*)\,. 
\end{equation}
So if we further assume that the initial perturbation $v_0 \in 
X$ is small enough so that 
\[
  2 C_1 \|v_0\|_X \,<\, 1\,,\quad \hbox{and} \quad
  4 C_1 C_2 K(1)\|v_0\|_X \,<\, 1\,,
\]
then it follows from \eqref{e:phibd} that $\phi(t) \le 2 C_1 
\|v_0\|_X < 1$ for all $t \in [0,T_*)$. Since $[0,T_*)$ is the maximal 
existence interval, this bound implies that $T_* = +\infty$ and that 
the solution of \eqref{e:pertv} satisfies
\[
  \sup_{t\ge 0}\|v(t)\|_{L^1} + \sup_{t \ge 0}(1+t)^{n/2}
  \|v(t)\|_{L^\infty} \,\le\,  2C_1 \|v_0\|_X\,.   
\]
This concludes the proof of Theorem~\ref{t:1} in the 
high-dimensional case $n \ge 3$. \hfill$\Box$ 

\section{Reduction to a normal form}\label{s:nf}

In low space dimensions the argument presented in the previous section
fails, and we must therefore have a closer look at the structure of
the perturbation equation. The idea is to introduce a normal form
transformation which simplifies the ODE dynamics in a neighborhood
of the periodic orbit $u_*$. Applying this transformation to the
reaction-diffusion equation \eqref{rdsystem}, we obtain a quasilinear
parabolic system which will be the starting point of our stability
analysis in Sections~\ref{s:le} and \ref{s:nl}.

We thus consider the ordinary differential equation
\begin{equation}\label{e:ode}
  u_t \,=\, f(u)\,, \quad u\in\R^N\,,
\end{equation}
with smooth nonlinearity $f\in C^\infty(\R^N,\R^N)$, and we assume the
existence of a time-periodic solution $u_*(t)=u_*(t+T)$ with minimal
period $T = 2\pi/\omega > 0$. As in Hypothesis~\ref{h:2}, we suppose
that $u_*$ is linearly asymptotically stable, in the sense that the
Floquet exponents $\lambda_1,\dots,\lambda_N$ are all contained in 
the open left half-plane, except for $\lambda_1 = 0$ (which is 
therefore algebraically simple). We shall show that the dynamics 
of \eqref{e:ode} near the periodic orbit $u_*$ is conjugate to the 
dynamics of the following normal form
\begin{equation}\label{e:odenf}
  \theta_t \,=\, \omega\,,\quad \tilde v_t \,=\, g(\theta,\tilde v)\,,
  \quad \theta\in S^1 \cong \R/2\pi\Z\,,\quad \tilde v\in B_\epsilon
  \subset \R^{N-1}\,,
\end{equation}
where $B_\epsilon$ denotes the open ball of radius $\epsilon > 0$ 
centered at the origin in $\R^{N-1}$. Here the vector field $g$ has 
the expansion
\begin{equation}\label{e:gnfdef}
  g(\theta,\tilde v) \,=\, L(\theta)\tilde v + g_2(\theta,\tilde v)
  [\tilde v,\tilde v]\,,
\end{equation}
where $L(\theta)$ is a real $(N-1)\times (N-1)$ matrix depending
smoothly on $\theta$, and $g_2(\theta,\tilde v)$ is a symmetric
bilinear form on $\R^{N-1}$ depending smoothly on $\theta,\tilde v$.
In particular $g(\theta,0) = 0$, hence \eqref{e:odenf} has a trivial
solution $\theta(t) = \omega t$, $\tilde v(t) = 0$ which will correspond
to the periodic solution $u_*(t)$ of \eqref{e:ode}. By construction,
the Floquet exponents $\lambda_2,\dots,\lambda_N$ of the time-periodic 
linear operator $L(\omega t)$ are all contained in the open left 
half-plane.

In what follows we denote by $\Phi(t)$ the flow of \eqref{e:ode} in a 
neighborhood of $u_*$, and by $\Phi_\mathrm{nf}(t)$ the flow of 
\eqref{e:odenf} in a neighborhood of $S^1 \times \{0\}$. These
local flows are defined at least for $t \ge 0$. 

\begin{Lemma}[Normal form]\label{p:fol}
Assume that the periodic solution $u_*$ is linearly asymptotically 
stable. Then there exist $\epsilon > \epsilon' > 0$ and a smooth 
diffeomorphism $\Psi$ from the solid torus $S^1\times B_\epsilon$ to 
a tubular neighborhood of the periodic orbit $u_*$ such that the local
flow in $S^1\times \R^{N-1}$ defined on $S^1 \times B_{\epsilon'}$ by
\[
  \Phi_\mathrm{nf}(t) \,=\, \Psi^{-1} \circ \Phi(t) \circ \Psi\,,
  \quad t \ge 0\,,
\]
is the flow induced by an ODE of the form \eqref{e:odenf}, 
\eqref{e:gnfdef}.
\end{Lemma}

\begin{Proof}
Since the periodic orbit $u_*$ is linearly asymptotically stable, we
can find a tubular neighborhood which is smoothly foliated by strong
stable fibers. Straightening out these fibers gives the desired
representation of the flow. For completeness we construct this
straightening change of coordinates in detail.

We start with the linearized equation at the periodic orbit,
$u_t=f'(u_*(t))u$, which possesses a linear invariant smooth
foliation: if we parametrize the periodic orbit $u_*(t)$ using $\theta
= \omega t \in S^1$, Floquet theory gives smooth families of
complementary subspaces $E^\mathrm{ss}(\theta)$ and $E^\mathrm{c}
(\theta)$, such that $\mathrm{dim}\,E^\mathrm{ss}(\theta)=N-1$ and 
$E^\mathrm{c}(\theta) = \mathrm{span}\,(u_*'(\theta/\omega))$. 
The linearized evolution leaves these subspaces invariant: 
$u(t)\in E^\mathrm{ss}(\theta)$ implies $u(t+\tau)\in 
E^\mathrm{ss}(\theta+\omega\tau)$, and the same holds for 
$E^\mathrm{c}(\theta)$. In particular, the family $\{E^\mathrm{ss}
(\theta)\}_{\theta \in S^1}$ forms a smooth normal bundle to the
periodic orbit $u_*$ (which is an orientable manifold in ambient 
Euclidean space), and such a bundle is necessarily trivial. 
Thus, we can find smooth coordinates $(\theta,v)\in S^1\times
\R^{N-1}$ and a smooth map $\Psi_0 : S^1\times\R^{N-1} \to \R^N$ 
such that $\Psi_0(\theta,0) = u_*(\theta/\omega)$ and 
$\Psi(\theta,\R^{N-1})=u_*(\theta/\omega)+E^\mathrm{ss}
(\theta)$ for all $\theta \in S^1$. 

On the other hand, for each $\theta \in S^1$, the strong stable 
manifold $W^\mathrm{ss}(\theta)$ of the nonlinear system 
\eqref{e:ode} is the graph of a local map $h_\theta : E^\mathrm{ss}
(\theta) \to E^\mathrm{c}(\theta)$, with $h_\theta(0)=0$ and 
$h_\theta'(0)=0$. The strong stable manifolds depend smoothly on the base point. In other words, $h_\theta$ depends smoothly on 
$\theta$, so that the map
\[
  u_*(\theta/\omega)+v^\mathrm{ss} ~\mapsto~
  \Psi_1(u_*(\theta/\omega)+v^\mathrm{ss}) \,:=\, 
  u_*(\theta/\omega) + v^\mathrm{ss}+h_\theta(v^\mathrm{ss})\,,
\]
defines a smooth diffeomorphism in a tubular neighborhood 
$\mathcal{U}$ of the periodic solution. Thus, if $\epsilon > 0$ 
is sufficiently small, the map $\Psi := \Psi_1 \circ \Psi_0 : 
S^1 \times B_\epsilon \to \mathcal{U}$ is also a smooth diffeomorphism
onto its image, and $\Psi(\theta,B_\epsilon) \subset W^\mathrm{ss}
(\theta)$ for all $\theta \in S^1$. Since
\begin{equation}\label{e:fol}
  \Phi(t )(W^\mathrm{ss}(\theta))\cap \mathcal{U} \,\subset\, 
  W^\mathrm{ss}(\theta+\omega t)\,, 
\end{equation}
we deduce that, if $\theta \in S^1$ and $\tilde v \in B_{\epsilon'}$ 
for some small $\epsilon'$, then $(\Phi(t)\circ\Psi)(\theta,\tilde v)$
belongs to the image of $\Psi$ for all $t \ge 0$, and 
\[
  \Phi_\mathrm{nf}(t)(\theta,\tilde v)  \,:=\, (\Psi^{-1}\circ\Phi(t)
  \circ\Psi)(\theta,\tilde v) \,=\, (\theta+\omega t,\hat v)\,,
\]
for some $\hat v \in B_\epsilon$. This immediately implies the
trivial form $\theta_t = \omega$ for the evolution equation associated to
$\Phi_\mathrm{nf}$. Moreover, by construction, $\Phi_\mathrm{nf}
(\theta,0) = (\theta+\omega t,0)$ for all $t \ge 0$, hence the
transverse variable $\tilde v$ evolves according to an ODE of the form
\eqref{e:odenf}, where the vector field satisfies $g(\theta,0) = 0$
and can therefore be expanded as in \eqref{e:gnfdef}.
\end{Proof}

We conclude this section with the transformation of the full
reaction-diffusion system. The pointwise change of coordinates
$u=\Psi(v)=\Psi(\theta,\tilde v)$ yields
\[
  v_t \,=\, \Psi'(v)^{-1}D\Delta \left(\Psi(v)\right)
  +f_\mathrm{nf}(v)\,, \quad f_\mathrm{nf}(v) \,=\, (\omega, g(v))^T\,,
\]
which can be expanded into 
\begin{equation}\label{e:ednf}
  v_t \,=\, \Psi'(v)^{-1}D\Psi'(v) \Delta v + \Psi'(v)^{-1} D\Psi''(v)
  [\nabla v,\nabla v] + f_\mathrm{nf}(v)\,.
\end{equation}
We are interested in the stability of the periodic orbit
$v_*(t)=(\omega t,0)^T$, and therefore we set $v = v_*(t)+w(t,x)$, so
that $w$ solves
\begin{equation}\label{e:ednfp}
  w_t \,=\, \Psi'(v_*+w)^{-1}D\Psi'(v_*+w) \Delta w + 
  \Psi'(v_*+w)^{-1}D\Psi''(v_*+w)[\nabla w,\nabla w] 
  +f_\mathrm{nf}^0(v_*+w)\,, 
\end{equation}
where now $f_\mathrm{nf}^0(v)=(0,g(v))^T$. In what follows, the first
component of the vector $w$ will play a distinguished role, as is 
clear from the expression of $f_\mathrm{nf}^0$. Thus we shall often 
write $w=(w_0,w_\mathrm{h})^T$, with $w_0 \in \R$ and $w_\mathrm{h}
\in \R^{N-1}$. 

\section{Linear evolution estimates}\label{s:le}

We consider the linearization of \eqref{e:ednfp} at $w = 0$, 
which reads
\begin{equation}\label{e:ednfl}
  w_t \,=\, \Psi'(v_*)^{-1}D\Psi'(v_*) \Delta w 
  +f_\mathrm{nf}'(v_*)w\,.
\end{equation}
To simplify the notations, we define 
\[
  A(t) \,=\, \Psi'(v_*(t))~, \quad \hbox{and} \quad 
  B(t) \,=\, \begin{pmatrix}0 & 0 \\ 0 & L(\omega t)
  \end{pmatrix}~,
\]
where $L(\theta)$ is the $(N-1) \times (N-1)$ matrix which 
appears in \eqref{e:gnfdef}. Note that $A(t)$, $B(t)$ are $T$-periodic 
$N \times N$ matrices, and that $A(t)$ is invertible for all $t$. The
linearization \eqref{e:ednfl} then becomes
\begin{equation}\label{e:ednfls}
  w_t \,=\, A(t)^{-1}DA(t) \Delta w + B(t)w\,, \quad t \in \R\,.
\end{equation}
By Fourier duality this system is equivalent to the family of 
ODEs
\begin{equation}\label{e:lf}
  w_t \,=\, -k^2 A(t)^{-1}DA(t)w + B(t)w\,, \quad t \in \R\,,
  \quad k \in \R^n\,.
\end{equation}
Since \eqref{e:lf} is related to \eqref{e:linf} by the $T$-periodic 
linear transformation $u=\Psi'(v_*(t))w$, it is clear that the 
Floquet spectrum of \eqref{e:lf} is identical to that of 
\eqref{e:linf} and therefore satisfies Hypothesis~\ref{h:2}.  
Let $M(t,s;k)$ denote the evolution operator defined by 
\eqref{e:lf}, so that any solution of \eqref{e:lf} satisfies 
$w(t) = M(t,s;k)w(s)$ for $t \ge s$. As $w=(w_0,w_\mathrm{h})^T \in 
\R\times\R^{N-1}$, it is natural to decompose the matrix $M$ 
in blocks as follows: 
\[
  M(t,s;k) \,=\, \begin{pmatrix} M_{00}(t,s;k) & M_{0\mathrm{h}}(t,s;k)
  \\ M_{\mathrm{h}0}(t,s;k) & M_{\mathrm{hh}}(t,s;k) \end{pmatrix}\,,
\]
where $M_{00}$, $M_{0\mathrm{h}}$, $M_{\mathrm{h}0}$,
$M_{\mathrm{hh}}$ are matrices of size $1\times 1$, $1\times (N-1)$, 
$(N-1)\times 1$, $(N-1)\times (N-1)$, respectively. The main 
result of this section is the following pointwise estimate 
on $M(t,s;k)$:

\begin{Proposition}[Pointwise estimates]\label{p:1}
There exist constants $C,d>0$ such that, for all $t \ge s$
and all $k \in \R^n$, one has
\begin{align}
  |M_{00}(t,s;k)| &\le C\,\rme^{-dk^2(t-s)}\,,\label{e:e1}\\
  |M_{0\mathrm{h}}(t,s;k)|+|M_{\mathrm{h}0}(t,s;k)| &\le
  \frac{C}{1+t-s} \,\rme^{-dk^2(t-s)}\,,\label{e:e2}\\
  |M_{\mathrm{hh}}(t,s;k)|&\le \frac{C}{(1+t-s)^2} \,\rme^{-dk^2(t-s)}\,,
  \label{e:e3}
\end{align}
where the norms on the left-hand side are arbitrary, $k$-independent 
matrix norms. 
\end{Proposition}

\begin{Proof}
Since the coefficients in \eqref{e:lf} are $T$-periodic, we have 
$M(t+T,s+T;k) = M(t,s;k)$ for all $t,s \in \R$ and all $k \in \R^n$. 
As a consequence, if $t \ge s$ and if $\tau_1,\tau_2 \in [0,T)$ 
are such that $t-\tau_1$ and $s+\tau_2$ are integer multiples 
of $T$, we have the identity
\begin{equation}\label{e:factor}
  M(t,s;k) \,=\, M(t,t-\tau_1;k)\,M(k)^m\,M(s+\tau_2,s;k)\,,
  \quad k \in \R^n\,,
\end{equation}
where $M(k) = M(T,0;k)$ and $m \in \N$ is such that $t-s = 
\tau_1 + mT + \tau_2$. To prove Proposition~\ref{p:1}, it is 
therefore sufficient to estimate $M(k)^m$ and $M(t,s;k)$ for
$0 \le t-s \le T$. 

\noindent
\textbf{Step 1:} Estimates on $M(k)^m$. \\
The general strategy is to distinguish between various parameters 
regimes. For small and intermediate values of $k$, we essentially
exploit Hypothesis~\ref{h:2}, while for large $k$ it is sufficient
to use the parabolicity of \eqref{e:lf}. 

\noindent
\textit{Small $k$:} We solve the ODE \eqref{e:lf} perturbatively
for $t \in [0,T]$ and obtain
\[
  w(t) \,=\, U(t,0)w(0) - k^2 \int_0^t U(t,s)A(s)^{-1} D A(s)
  U(s,0)w(0)\,\rmd s + \rmO(k^4)\,,
\]
where $U(t,s) = M(t,s;0)$ is the evolution operator associated 
to the equation $w_t = B(t)w$. In particular, setting $t = T$, 
we find
\begin{align}\nonumber
  M(k) \,&=\, U(T,0)\left(\id - k^2 \int_0^T (A(t)U(t,0))^{-1}
  D (A(t)U(t,0))\,\rmd t + \rmO(k^4)\right) \\ \label{e:Mkexp}
  \,&=\, \begin{pmatrix}1 & 0 \\ 0 & V\end{pmatrix}
  \begin{pmatrix}1 -d_0 T k^2 + \rmO(k^4) & \rmO(k^2)\\
  \rmO(k^2) & \id + \rmO(k^2)\end{pmatrix}\,,
\end{align}
where $V$ is the $(N-1)\times(N-1)$ Floquet matrix associated
to the $T$-periodic linear operator $L(\omega t)$, and 
\begin{equation}\label{e:d0def2}
  d_0 \,=\, \frac1T \int_0^T e_1^T (A(t)U(t,0))^{-1}
  D (A(t)U(t,0))e_1\,\rmd t \,=\,\frac1T \int_0^T e_1^T A(t)^{-1}
  D A(t)e_1\,\rmd t\,.
\end{equation}
Here $e_1 = (1,0)^T$ is the first vector of the canonical basis in
$\R^N$. As was already observed, all eigenvalues of $V$ are contained
in the disk $\{z \in \C\,|\, |z| < e^{-\nu}\}$ for some $\nu > 0$, and
it follows from \eqref{e:Mkexp} that $M(k)$ has exactly $N-1$
eigenvalues in this disk if $k$ is sufficiently small. The remaining
Floquet multiplier has the expansion $1 - d_0 T k^2 + \rmO(k^4)$, in
agreement with Hypothesis~\ref{h:2}. Incidentally, we observe that
\eqref{e:d0def2} is identical to \eqref{e:d0def}. Indeed, since
$u_*(t) = \Psi(v_*(t)) = \Psi(\omega t e_1)$, we have $u_*'(t) =
\omega\Psi'(v_*(t))e_1 = \omega A(t)e_1$, and it is also
straightforward to verify that the bounded solution of the adjoint
equation \eqref{e:adj} is $U_*(t) = (A(t)^{-1})^T e_1$. Thus
\eqref{e:d0def2} can be written as
\[
  d_0 \,=\, \frac{1}{\omega T} \int_0^T U_*(t)^T D u_*'(t)\,\rmd t
  \,=\, \frac{\int_0^T U_*(t)^T D u_*'(t)\,\rmd t}
  {\int_0^T U_*(t)^T u_*'(t)\,\rmd t}\,,
\]
and Hypothesis~\ref{h:2} guarantees that $d_0 > 0$. 

For $k \in \R^n$ sufficiently small, let $P(k)$ denote the spectral
projection onto the one-dimensional eigenspace of $M(k)$ corresponding
to the neutral Floquet exponent $\lambda_1(k) = -d_0 k^2 + \rmO(k^4)$.
From \eqref{e:Mkexp} it is easy to verify that $P(k)$ has the 
following form
\[
  P(k) \,=\, \frac{1}{1+k^4 b^T a}\begin{pmatrix} 1 & k^2 b^T \\
  k^2 a & k^4 a b^T\end{pmatrix}\,,
\]
where $a(k), b(k)$ are $(N-1)$-dimensional vectors with $a(k), 
b(k) = \rmO(1)$ as $k \to 0$. By construction, we have for 
any $m \in \N^*$:
\begin{equation}\label{e:Mkm}
  M(k)^m \,=\, M(k)^m P(k) + M(k)^m (\id - P(k)) \,=\,
  \rme^{mT\lambda_1(k)}P(k) + (M(k) (\id - P(k)))^m\,. 
\end{equation}
Since $\lambda_1(k) \le -d_1 k^2$ for small $k$ if $0 < d_1 < d_0$, 
and since the spectral radius of $M(k)(\id - P(k))$ is smaller 
than $e^{-\nu}$, we conclude that
\[
  |M(k)^m| \,\le\, C\,\rme^{-d_1 k^2 mT}\begin{pmatrix} 1 & k^2 \\
  k^2 & k^4\end{pmatrix} + \rmO(e^{-\nu m}) \,\le\, 
  C\,\rme^{-d_1 k^2 mT}\begin{pmatrix} 1 & (mT)^{-1} \\
  (mT)^{-1} & (mT)^{-2}\end{pmatrix}\,,
\]
for all $m \in \N^*$ if $|k| \le \kappa_0 \ll 1$. Here the 
matrix norm $|\cdot|$ is applied separately to each of the 
four blocks of $M(k)^m$, and for convenience the four upper
bounds are collected in a $2\times 2$ matrix. 

\noindent
\textit{Large $k$:} In this parameter regime, it is more 
convenient to set $v(t) = A(t)w(t)$ and to solve the $v$-equation
corresponding to \eqref{e:lf}, namely
\begin{equation}\label{e:lff}
  v_t \,=\, -k^2 D v + C(t)v\,, \quad \hbox{where} \quad
  C(t) \,=\, A'(t)A(t)^{-1} + A(t)B(t)A(t)^{-1}\,.
\end{equation}
The matrix $C(t)$ is $T$-periodic, hence uniformly bounded. 
A standard energy estimate yields
\[
  \frac12\frac{\rmd}{\rmd t}\|v(t)\|_2^2 \,=\, 
  -k^2 (v(t),Dv(t)) + (v(t),C(t)v(t)) \,\le\, -d_2 k^2 
  \|v(t)\|_2^2 + K \|v(t)\|_2^2\,,
\]
where $d_2 > 0$ is the smallest eigenvalue of the (symmetric 
and positive) matrix $D$, and $K = \sup_{t \in [0,T]}\|C(t)\|_2$. 
Thus any solution of \eqref{e:lff} satisfies $\|v(t)\|_2 \le 
\rme^{(-d_2 k^2 + K)t}\|v(0)\|_2$, and returning to the $w$-equation
we obtain
\begin{equation}\label{e:largew}
  \|w(t)\|_2 \,\le\, C\,\rme^{(-d_2 k^2 + K)t}\,\|w(0)\|_2\,,
  \quad t \ge 0\,,
\end{equation}
for some $C > 0$. In particular, if we choose $t = mT$ and if 
we assume that $|k| \ge \kappa_1$ with $\kappa_1 = (2K/d_2)^{1/2}$, 
we arrive at
\[
  |M(k)^m| \,\le\, C\,\rme^{-d_3 k^2 mT}\,, \quad \hbox{where}
  \quad d_3 \,=\, \frac{d_2}{2T}\,.
\] 

\noindent
\textit{Intermediate $k$:} By Hypothesis~\ref{h:2}, if $\kappa_0 
\le |k| \le \kappa_1$, the spectrum of $M(k)$ is entirely contained
in the disk $\{z \in \C\,|\, |z| < \rme^{-\mu}\}$ for some 
$\mu > 0$. Moreover, the resolvent matrix $(z - M(k))^{-1}$ is 
uniformly bounded for all $z$ on the circle $\{|z| = \rme^{-\mu}\}$
and all $k$ in the annulus $\kappa_0 \le |k| \le \kappa_1$. If 
$0 < d_4 < \mu/(\kappa_1^2T)$, it follows that
\[
  |M(k)^m| \,\le\, C\,\rme^{-d_4 k^2 mT}\,, \quad m \in \N\,,
\]
where the constant $C$ is independent of $k$. 

Summarizing the results obtained so far, we proved that there
exist $C > 0$ and $d > 0$ such that
\begin{equation}\label{e:step1}
  |M(k)^m| \,\le\, C\,\rme^{-d k^2 mT} \begin{pmatrix} 1 & (mT)^{-1} \\
  (mT)^{-1} & (mT)^{-2}\end{pmatrix}\,,
\end{equation}
for all $k \in \R^n$ and all $m \in \N^*$. This is the particular 
case of \eqref{e:e1}--\eqref{e:e3} when $s = 0$ and $t = mT$. 

\noindent
\textbf{Step 2:} Estimates on $M(t,s;k)$ for $0 \le t-s \le 2T$. \\
Our goal is to show that
\begin{equation}\label{e:step2}
  |M(t,s;k)| \,\le\, C\,\rme^{-d k^2(t-s)} \begin{pmatrix} 1 & 
  k^2(t-s) \\ k^2(t-s) & 1\end{pmatrix}\,,
\end{equation}
for some $C > 0$ and $d > 0$. Note that \eqref{e:step2} implies
\eqref{e:e1}--\eqref{e:e3} if $0 \le t-s \le 2T$. 

In the case where $k^2(t-s)$ is small, say $k^2(t-s) \le \kappa_2 
\ll 1$, we can solve \eqref{e:lf} perturbatively and obtain
as in \eqref{e:Mkexp}
\begin{equation}\label{e:Mtsexp}
  M(t,s;k) \,=\, \begin{pmatrix} 1 & 0 \\ 0 & V(t,s)\end{pmatrix}
  \,\Bigl(\id + \rmO(k^2(t-s))\Bigr)\,,
\end{equation}
from which \eqref{e:step2} follows (for any fixed $d > 0$). If 
$k^2(t-s) \ge \kappa_2$, then using \eqref{e:largew} we immediately 
find
\[
  |M(t,s;k)| \,\le\, C\,\rme^{(-d_2 k^2 + K)(t-s)} \,\le\, 
  C\,\rme^{2KT}\,\rme^{-d_2 k^2 (t-s)}\,,
\]
which implies \eqref{e:step2} with $d = d_2$. 

It is now straightforward to conclude the proof of 
Proposition~\ref{p:1}. In view of \eqref{e:step2}, it remains 
to prove \eqref{e:e1}--\eqref{e:e3} for $t-s \ge 2T$. 
We decompose $t-s = \tau_1 + mT + \tau_2$ with $\tau_1,\tau_2 
\in [0,T)$ and $m \in \N^*$, and we factorize $M(t,s;k)$ 
as in \eqref{e:factor}. Using \eqref{e:step1}, \eqref{e:step2}, 
we find
\begin{align*}
  |M(t,s;k)| \,&\le\, C\,\rme^{-dk^2 (t-s)}
  \begin{pmatrix} 1 & k^2 \tau_1 \\ k^2 \tau_1 & 1\end{pmatrix}
  \begin{pmatrix} 1 & (mT)^{-1} \\ (mT)^{-1} & (mT)^{-2}\end{pmatrix}
  \begin{pmatrix} 1 & k^2 \tau_2 \\ k^2 \tau_2 & 1\end{pmatrix} \\
  \,&\le\, C\,\rme^{-dk^2 (t-s)} \begin{pmatrix} 1 & 
  (1+t-s)^{-1} \\ (1+t-s)^{-1} & (1+t-s)^{-2}\end{pmatrix}\,,
\end{align*}
which is the desired estimate. 
\end{Proof}

\begin{Remark}\label{r:asym}
For all $k \in \R^n$ and all $s \in \R$, we have
\[
  \lim_{t \to +\infty} M(t,s;kt^{-1/2}) \,=\, \rme^{-d_0 k^2}
  \begin{pmatrix} 1 & 0 \\ 0 & 0\end{pmatrix}\,,
\]
where $d_0$ is as in Hypothesis~\ref{h:2}. This follows from the
proof of Proposition~\ref{p:1}, and in particular from 
\eqref{e:factor}, \eqref{e:Mkm}, \eqref{e:Mtsexp}. 
\end{Remark}

\begin{Proposition}[Estimates for the derivatives]\label{p:2}
For any $\alpha \in \N^n$ there exists $C > 0$ such that, for all 
$t > s$ and all $k \in \R^n$, 
\begin{equation}\label{e:derest}
  |\partial_k^\alpha M(t,s;k)| \,\le\, C (t-s)^{|\alpha|/2}
  \,\rme^{-dk^2(t-s)}\,.
\end{equation}
The same estimate holds for $M_{00}$, $(1+t-s)M_{0\mathrm{h}}$, 
$(1+t-s)M_{\mathrm{h}0}$, and $(1+t-s)^2 M_{\mathrm{hh}}$.
\end{Proposition}

\begin{Proof}
It is clear from \eqref{e:lf} that $M(t,s;k)$ depends on the 
parameter $k \in \R^n$ only through the scalar quantity $p = k^2$. 
In this proof, we set $M(t,s;k) = \tilde M(t,s;k^2)$ and we
consider the derivatives of $\tilde M(t,s;p)$ with respect to $p$. 
Our goal is to prove the estimate
\begin{equation}\label{e:derest2}
  |\partial_p^j \tilde M(t,s;p)| \,\le\, C_j (t-s)^j\,\rme^{-dp(t-s)}, 
  \quad j \in \N\,,
\end{equation}
which implies immediately \eqref{e:derest}. Differentiating
\eqref{e:lf} with respect to $p = k^2$, we find
\[
  (\partial_p w)_t \,=\, -\mathcal{D}(t)  w + (B(t) - p\mathcal{D}(t))
  (\partial_p w)\,, 
\]
where $\mathcal{D}(t) = A(t)^{-1}DA(t)$. Since $\partial_p 
M(s,s;p) = 0$, we deduce that
\[
  \partial_p\tilde M(t,s;p) \,=\, -\int_s^t \tilde M(t,t_1;p)
  \mathcal{D}(t_1)\tilde M(t_1,s;p)\,\rmd t_1\,. 
\] 
Iterating this procedure, we obtain for any $j \in \N$ the 
representation formula
\begin{align*}
 \partial_p^j\tilde M(t,s;p) &= (-1)^j j! \int_s^t \int_s^{t_1}
 \dots \int_s^{t_{j-1}} \tilde M(t,t_1;p)\mathcal{D}(t_1) \times\\
 &\quad \times M(t_1,t_2;p)\mathcal{D}(t_2) \dots 
 M(t_{j-1},t_j;p)\mathcal{D}(t_j)M(t_j,s;p)\,\rmd t_j \dots \,\rmd
 t_1\,. 
\end{align*}
Using now the pointwise estimates established in Proposition~\ref{p:1}, 
we easily obtain \eqref{e:derest2}. Similar bounds can be proved
for $\tilde M_{00}$, $(1+t-s)\tilde M_{0\mathrm{h}}$, $(1+t-s)\tilde
M_{\mathrm{h}0}$, and $(1+t-s)^2 \tilde M_{\mathrm{hh}}$ (we omit the
details). 
\end{Proof}

We now convert our estimates on the Fourier multipliers $M(t,s;k)$ 
into bounds on the linear evolution equation \eqref{e:ednfls}
in various $L^p$ spaces. The two-parameter semigroup 
$\mathcal{M}(t,s)$ associated to \eqref{e:ednfls} is defined 
using Fourier transform by the relation
\begin{equation}\label{e:le}
  \widehat{(\mathcal{M}(t,s)v)}(k) \,=\, M(t,s;k)\hat{v}(k)\,,
  \quad k \in \R^n\,.
\end{equation}
The following proposition contains the main estimates on 
$\mathcal{M}(t,s)$ which will be used in the nonlinear stability
proof. 

\begin{Proposition}[$L^p$-$L^q$ estimates]\label{p:3}
There exist a positive constant $C$ such that, for all $t > s$ and 
all $1 \le p \le q \le \infty$, one has
\begin{equation}\label{e:lplq}
  \|\mathcal{M}(t,s)v\|_{L^q(\R^n)} \,\le\, \frac{C}{(t-s)^{\frac{n}2
  (\frac1p-\frac1q)}}\ \|v\|_{L^p(\R^n)}\,.
\end{equation}
The same estimate holds for $\mathcal{M}_{00}$, $(1+t-s)
\mathcal{M}_{0\mathrm{h}}$, $(1+t-s)\mathcal{M}_{\mathrm{h}0}$, and 
$(1+t-s)^2\mathcal{M}_{\mathrm{hh}}$.
\end{Proposition}

\begin{Proof}
By construction $\mathcal{M}(t,s)$ is the convolution operator 
with the function $x \mapsto \mathcal{M}(t,s;x)$, which is 
just the inverse Fourier transform of $k \mapsto M(t,s;k)$. 
Thus using the pointwise estimate \eqref{e:e1}, we easily obtain
\[
  \|\mathcal{M}(t,s;\cdot)\|_{L^\infty(\R^n)} \,\le\, C 
  \|M(t,s;\cdot)\|_{L^1(\R^n)} \le \frac{C}{(t-s)^{n/2}}\,.
\]
To estimate the $L^1$ norm of $\mathcal{M}(t,s;\cdot)$, we use 
Sobolev embeddings. Let $m \in \N$ be the smallest integer 
such that $m > n/2$. Using the estimates of Proposition~\ref{p:2}
together with H\"older's inequality and Parseval's identity, 
we find
\begin{align*}
  &\|\mathcal{M}(t,s;\cdot)\|_{L^1(\R^n)} \,=\, 
  \int_{\R^n} \Bigl(1 + \frac{|x|^2}{t-s}\Bigr)^{-m/2}
  \Bigl(1 + \frac{|x|^2}{t-s}\Bigr)^{m/2} |\mathcal{M}(t,s;x)|\ 
  \rmd x\\
  &\qquad \,\le\, C(t-s)^{n/4}\left(\int_{\R^n}\Bigl(1 + \frac{|x|^2}{t-s} 
  \Bigr)^m |\mathcal{M}(t,s;x)|^2\,\rmd x\right)^{1/2}\\
  &\qquad \,\le\, C(t-s)^{n/4} \left(\int_{\R^n}\sum_{|\alpha| \le m} 
  (t-s)^{-|\alpha|}|\partial_k^\alpha M(t,s;k)|^2 \rmd k\right)^{1/2} 
  \le C\,.
\end{align*}
Summarizing, we have shown that 
\[
  \|\mathcal{M}(t,s;\cdot)\|_{L^p(\R^n)} \,\le\, \frac{C}{(t-s)^{
  \frac{n}2(1-\frac1p)}}, \quad 1 \le p \le \infty\,,
\]
and \eqref{e:lplq} follows by Young's inequality. The estimates
for $\mathcal{M}_{00}$, $(1+t-s)\mathcal{M}_{0\mathrm{h}}$, 
$(1+t-s)\mathcal{M}_{\mathrm{h}0}$, and $(1+t-s)^2\mathcal{M}_{
\mathrm{hh}}$ are proved in exactly the same way. 
\end{Proof}

\begin{Remark}\label{r:1}
Under the assumptions of Proposition~\ref{p:3}, we also have
\[
  \|\partial_x^\alpha \mathcal{M}(t,s)v\|_{L^q(\R^n)} \,\le\, 
  \frac{C_\alpha}{(t-s)^{\frac{n}2(\frac1p-\frac1q)+\frac{|\alpha|}{2}}}\,
  \|v\|_{L^p(\R^n)}\,, \quad \alpha \in \N^n\,,
\]
and the same estimates hold for $\mathcal{M}_{00}$, $(1+t-s)
\mathcal{M}_{0\mathrm{h}}$, $(1+t-s)\mathcal{M}_{\mathrm{h}0}$, and 
$(1+t-s)^2\mathcal{M}_{\mathrm{hh}}$. This is obvious in view 
of Proposition~\ref{p:1}, since the operator $\partial_x^\alpha 
\mathcal{M}(t,s)$ is the Fourier multiplier associated to 
the function $(ik)^\alpha M(t,s;k)$. 
\end{Remark} 

\begin{Remark}\label{r:2}
For later use, we also observe that, if $t > 0$ and $0 < s < t$, 
then
\begin{equation}\label{e:lplqdiff}
  \|(\mathcal{M}_{00}(t,s) - \mathcal{M}_{00}(t,0))v\|_{L^q(\R^n)} 
  \,\le\, \frac{C}{(t-s)^{\frac{n}2(\frac1p-\frac1q)}}\,\frac{s}{t}\, 
  \|v\|_{L^p(\R^n)}\,,
\end{equation}
for $1 \le p \le q \le \infty$. If $t/2 \le s \le t$, this bound
follows immediately from \eqref{e:lplq} and the triangle 
inequality. If $0 < s \le t/2$, we observe that
\[
  M_{00}(s,0;k) \,=\, 1 + k^2 s R_0(s,k)\,, \quad 
  M_{\mathrm{h}0}(s,0;k) \,=\, k^2 s R_\mathrm{h}(s,k)\,, 
\]
where $R_0(s,k)$, $R_\mathrm{h}(s,k)$ are uniformly bounded for
$s > 0$ and $k \in \R^n$, see \eqref{e:Mtsexp}. Since
\[
  M_{00}(t,s;k) - M_{00}(t,0;k) \,=\, M_{00}(t,s;k)(1 - M_{00}(s,0;k))
  - M_{0\mathrm{h}}(t,s;k)M_{\mathrm{h}0}(s,0;k)\,,
\]
we obtain the pointwise bound
\[
  |M_{00}(t,s;k) - M_{00}(t,0;k)| \,\le\, C k^2 s \,\rme^{-k^2 d(t-s)}
  \,\le\, C\,\frac{s}{t-s}\,\rme^{-k^2 d(t-s)}\,,
\]
which allows to establish the $L^p$-$L^q$ estimate \eqref{e:lplqdiff}
using the same arguments as in the proof of Proposition~\ref{p:3}. 
\end{Remark}

Finally, to control the quasilinear terms in the perturbation
equation we will use maximal regularity properties of the 
evolution semigroup $\mathcal{M}(t,s)$. 

\begin{Proposition}[Maximal regularity of type $L^r$]\label{p:4}
For any $r \in (1,+\infty)$ and any $\tilde T > 0$, there exists 
$C > 0$ such that the following holds. If $v \in L^r((0,\tilde T),
L^2(\R^n))$ and if $w$ satisfies
\begin{equation}\label{e:wvdef}
  w(t) \,=\, \int_0^t \mathcal{M}(t,s)v(s)\,\rmd s\,, \quad t \in 
  [0,\tilde T]\,,
\end{equation}
then
\begin{equation}\label{e:maxres}
  \int_0^{\tilde T} \|\Delta w(t)\|_{L^2}^r \,\rmd t \,\le\, C
  \int_0^{\tilde T} \|v(t)\|_{L^2}^r \,\rmd t\,.
\end{equation}
\end{Proposition}

\begin{Proof}
For any $t \in \R$ the generator $\mathcal{L}(t) = A(t)^{-1}DA(t) 
\Delta + B(t)$ is an elliptic operator in $L^2(\R^n)$ with 
(time-independent) domain $H^2(\R^n)$. Moreover $\mathcal{L}(t)$, 
considered as a bounded operator from $H^2(\R^n)$ into $L^2(\R^n)$,
is a smooth function of $t$. Thus estimate \eqref{e:maxres} is 
a particular case of the results established in \cite{HM00,PS01}. 
\end{Proof}

We also state a corollary which will be useful in the next section.

\begin{Corollary}\label{c:max}
Fix $r \in (1,\infty)$ and $\alpha \in \R$. There exists $C > 0$
such that, if 
\[
  w(t) \,=\, \int_{t-T}^t \mathcal{M}(t,s)v(s)\,\rmd s\,, \quad 
  t \ge T\,,
\]
then
\[
  \int_T^{\infty} (1+t)^\alpha \|\Delta w(t)\|_{L^2}^r \,\rmd t 
  \,\le\, C \int_0^\infty (1+t)^\alpha \|v(t)\|_{L^2}^r \,\rmd t\,.
\]
\end{Corollary}

\begin{Proof}
Fix $k \in \N^* := \N\setminus\{0\}$. If $t \in [kT,(k+1)T]$, 
where $T > 0$ is the period of $\mathcal{L}(t)$, we can write
\[
  w(t) \,=\, \int_{(k-1)T}^t \mathcal{M}(t,s)v(s)\,\rmd s
  - \int_{(k-1)T}^{t-T} \mathcal{M}(t,s)v(s)\,\rmd s\,.
\]
Since $\mathcal{M}(t+T,s+T) = \mathcal{M}(t,s)$, and since
$0 \le t - (k-1)T \le 2T$, we can use Proposition~\ref{p:4}
(with $\tilde T = 2T$) to control both terms in the 
right-hand side. We obtain
\[
  \int_{kT}^{(k+1)T} \|\Delta w(t)\|_{L^2}^r \,\rmd t \,\le\, C
  \int_{(k-1)T}^{(k+1)T} \|v(t)\|_{L^2}^r \,\rmd t\,,
\]
for some $C > 0$ independent of $k$. Multiplying both sides
by $(kT)^\alpha \approx (1+t)^\alpha$ and summing over 
$k \in \N^*$, we obtain the desired result. 
\end{Proof}

\section{Nonlinear stability in low dimensions}\label{s:nl}

This section is devoted to the proof of Theorem~\ref{t:1} in the case
$n \le 3$. Smoothness of the change of coordinates $\Psi$ implies
that it is sufficient to prove the decay estimate \eqref{e:t1decay} 
for the transformed equation \eqref{e:ednfp}, with smallness
assumptions on the initial data $w_0$. To simplify the notations, 
we rewrite \eqref{e:ednfp} in the compact form
\begin{equation}\label{e:13}
  w_t \,=\, \mathcal{L}(t) w + F(t,w,\Delta w) + G(t,w,\nabla w)+ 
  H(t,w)\,,
\end{equation}
where $\mathcal{L}(t)= A(t)^{-1}DA(t)\Delta + B(t)$ is the linear 
operator studied in Section~\ref{s:le} and 
\begin{align*}
  F(t,w,\Delta w)\,&=\,\left(\Psi'(v_*(t)+w)^{-1}-\Psi'(v_*(t))^{-1}\right)
  D \left(\Psi'(v_*(t)+w)\right)\Delta w \\
  &\quad + \Psi'(v_*(t))^{-1}D\left(\Psi'(v_*(t)+w)-\Psi'(v_*(t))\right)
  \Delta w\,,\\
  G(t,w,\nabla w) \,&=\, \Psi'(v_*(t)+w)^{-1}D\Psi''(v_*(t)+w)[\nabla w,
  \nabla w]\,,\\ 
  H(t,w) \,&=\, (0,\hat g_2(t,w))^T, \quad \hbox{where} \quad 
  \hat g_2(t,w) \,=\, g_2(v_*(t)+w)[w_\mathrm{h},w_\mathrm{h}]\,.
\end{align*}
Clearly, $F(t,w,\Delta w)=\rmO(w\Delta w)$ and $G(t,w,\nabla w)=
\rmO(|\nabla w|^2)$. More precisely, if $\|w\|_{L^\infty}$ is 
sufficiently small, we have
\[
\begin{array}{l}
  \|F(t,w,\Delta w)\|_{L^1} \,\le\,  C\|w\|_{L^2}\|\Delta w\|_{L^2}\,,\\[1mm]
  \|F(t,w,\Delta w)\|_{L^2} \,\le\,  C\|w\|_{L^\infty}\|\Delta w\|_{L^2}\,,
\end{array} \qquad
\begin{array}{l}
  \|G(t,w,\nabla w)\|_{L^1} \,\le\, C\|\nabla w\|^2_{L^2}\,,\\[1mm]
  \|G(t,w,\nabla w)\|_{L^2} \,\le\, C\|\nabla w\|^2_{L^4}\,.
\end{array}
\]
Since 
\[
  \|\nabla w\|^2_{L^2}\,\le\, C\|w\|_{L^2}\|\Delta w\|_{L^2}\quad 
  \hbox{and}\quad \|\nabla w\|^2_{L^4}\,\le\, C\|w\|_{L^\infty}
  \|\Delta w\|_{L^2}\,,
\]
we find
\begin{equation}\label{e:nlest}
\begin{array}{l}
  \|F(t,w,\Delta w)\|_{L^1}+\|G(t,w,\nabla w)\|_{L^1} \,\le\, 
  C\|w\|_{L^2}\|\Delta w\|_{L^2}\,,\\[1mm]
  \|F(t,w,\Delta w)\|_{L^2}+\|G(t,w,\nabla w)\|_{L^2} \,\le\, 
  C\|w\|_{L^\infty}\|\Delta w\|_{L^2}\,.
\end{array}
\end{equation}
Under the same assumptions, we also have
\begin{equation}
  \|H(t,w)\|_{L^1}\,\le\, C \|w_\mathrm{h}\|_{L^2}^2\quad \hbox{and}\quad
  \|H(t,w)\|_{L^2}\,\le\, C \|w_\mathrm{h}\|_{L^2}\|w_\mathrm{h}\|_{L^\infty}\,.
\label{e:nlesth}
\end{equation}
Setting $K(t,w,\nabla w,\Delta w)=F(t,w,\Delta w)+G(t,w,\nabla w)$,
we can write the integral equation associated with \eqref{e:ednfp} 
in the form,
\begin{equation}\label{e:mild}
  w(t)\,=\,\mathcal{M}(t,0)w_0 + \int_0^t \mathcal{M}(t,s)
  K(s,w,\nabla w,\Delta w) \,\rmd s + \int_0^t \mathcal{M}(t,s) 
  H(s,w) \,\rmd s\,. 
\end{equation}

We now describe the function space in which we shall look for
solutions of \eqref{e:mild}. Assume that $n \le 3$ and choose 
$r \in (4,+\infty)$, so that
\begin{equation}\label{e:rineq}
  \frac1r \,<\, \frac{n}4 \,<\, 1 - \frac1r\,.
\end{equation}
We define the Banach space
\[
  Y \,=\, \Bigl\{w \in C^0([0,+\infty),L^1(\R^n)\cap
  L^\infty(\R^n)) \cap L^r((0,+\infty),\dot H^2(\R^n)) \,\Big|\, 
  \|w\|_Y < \infty\Bigr\}\,,
\]
where
\begin{equation}\label{e:norm}
  \|w\|_Y \,=\, \sup_{t \ge 0} \|w(t)\|_{L^1} + \sup_{t \ge 0}
  (1+t)^{n/2} \|w(t)\|_{L^\infty} + \left(\int_0^\infty 
  (1+t)^r \|\Delta w(t)\|^r_{L^2}\,\rmd t\right)^{1/r}\,.
\end{equation}

\begin{Lemma}
If $w_0 \in L^1(\R^n) \cap H^2(\R^n)$, the linear solution $W_0(t) = 
\mathcal{M}(t,0)w_0$ belongs to $Y$ and $\|W_0\|_Y \le C_1
(\|w_0\|_{L^1} + \|w_0\|_{H^2})$ for some $C_1 > 0$. 
\end{Lemma}

\begin{Proof}
Since $n \le 3$, we have $H^2(\R^n) \hookrightarrow L^\infty(\R^n)$, 
hence $w_0 \in L^1(\R^n) \cap L^\infty(\R^n)$. Thus, using the linear 
estimates established in Proposition~\ref{p:3}, we obtain
\[
  \sup_{t \ge 0} \|W_0(t)\|_{L^1} + \sup_{t \ge 0}(1+t)^{n/2} 
  \|W_0(t)\|_{L^\infty} \,\le\, C(\|w_0\|_{L^1} + \|w_0\|_{L^\infty})\,.  
\]
On the other hand, using Proposition~\ref{p:3} and Remark~\ref{r:1}, 
we find
\[
  \|\Delta W_0(t)\|_{L^2} \,\le\, C \|w_0\|_{H^2} \quad \hbox{for } 
  t \le 1\,, \qquad 
  \|\Delta W_0(t)\|_{L^2} \,\le\, \frac{C \|w_0\|_{L^1}}{t^{1+n/4}}
  \quad \hbox{for } t \ge 1\,. \qquad 
\]
As $rn/4 > 1$ by \eqref{e:rineq}, we conclude that
\[
  \left(\int_0^\infty (1+t)^r \|\Delta W_0(t)\|^r_{L^2}\,\rmd 
  t\right)^{1/r} \,\le\, C(\|w_0\|_{L^1} + \|w_0\|_{H^2})\,.  
\]
\end{Proof}

The next step consists in estimating the integral terms in
\eqref{e:mild}, namely
\[
  \mathcal{I}(t) \,=\, \int_0^t \mathcal{M}(t,s)K(s,w,\nabla w,\Delta w) 
  \,\rmd s\,, \quad \hbox{and} \quad \mathcal{J}(t) \,=\, \int_0^t 
  \mathcal{M}(t,s) H(s,w) \,\rmd s\,.
\]

\begin{Proposition}\label{p:fp}
There exist $C_2 > 0$ and $\delta_2 > 0$ such that, for all $w \in Y$
with $\|w\|_Y \le \delta_2$, we have
\[
  \|\mathcal{I}\|_Y + \|\mathcal{J}\|_Y \,\le\, C_2\|w\|_Y^2\,.
\]
\end{Proposition}

\begin{Proof}
Throughout the proof $C$ denotes a constant that changes between 
estimates, but does not depend on $w$. Smallness of $w$ in $Y$ implies 
that estimates \eqref{e:nlest} and \eqref{e:nlesth} hold for the
nonlinearities. 

\noindent
\textbf{Estimate on $\|\mathcal{I}(t)\|_{L^1}$}\\
Using \eqref{e:lplq} with $p = q = 1$ and the first estimate in  
\eqref{e:nlest}, we find
\begin{align}\nonumber
  \|\mathcal{I}(t)\|_{L^1} &\le C\int_0^t \|w(s)\|_{L^2}\|\Delta 
  w(s)\|_{L^2}\,\rmd s\\ \label{e:IL1}
  &\le C \|w\|_Y\int_0^t \frac{1}{(1+s)^{\frac{n}{4}+1}}(1+s) 
  \|\Delta w(s)\|_{L^2}\,\rmd s\\ \nonumber
  &\le C \|w\|_Y \left(\int_0^t \frac{1}{(1+s)^{(\frac{n}{4}+1)
  \frac{r}{r-1}}}\,\rmd s\right)^{1-1/r} \left(\int_0^t (1+s)^r 
  \|\Delta w(s)\|_{L^2}^r \,\rmd s\right)^{1/r}\,.
\end{align}
In the second inequality we used the bound $\|w(s)\|_{L^2} \le 
\|w(s)\|_{L^1}^{1/2}\|w(s)\|_{L^\infty}^{1/2} \le \|w\|_Y
(1+s)^{-n/4}$, and in the last line H\"older's inequality. Taking 
the supremum over $t \ge 0$ and using \eqref{e:norm}, we conclude 
that 
\[
  \sup_{t \ge 0}\|\mathcal{I}(t)\|_{L^1} \,\le\, C \|w\|_Y^2\,.
\]

\noindent
\textbf{Estimate on $\|\mathcal{I}(t)\|_{L^\infty}$} \\
For $t\le 1$, we use \eqref{e:lplq} with $(p,q) = (2,\infty)$ 
and the second estimate in \eqref{e:nlest}. We obtain
\begin{equation}\label{e:Ismallt}
  \|\mathcal{I}(t)\|_{L^\infty} \le C\int_0^t
  \frac{1}{(t-s)^{n/4}}\|w(s)\|_{L^\infty}\|\Delta w(s)\|_{L^2}
  \,\rmd s \le C t^{1-\frac{n}4-\frac1r}\|w\|_Y^2\,, 
\end{equation}
where the last estimate is again a consequence of H\"older's 
inequality. Note that $1-\frac{n}4-\frac1r > 0$ by \eqref{e:rineq}.
For $t\ge 1$ we split
\begin{align*}
  \mathcal{I}(t) &= \int_0^{t/2}\mathcal{M}(t,s)
  K(s,w,\nabla w,\Delta w)\,\rmd s + 
  \int_{t/2}^t\mathcal{M}(t,s)K(s,w,\nabla w,\Delta w)
  \,\rmd s\\ 
  \,&=:\, \mathcal{I}_1(t) + \mathcal{I}_2(t)\,.
\end{align*}
Using \eqref{e:lplq} with $(p,q) = (1,\infty)$ and
proceeding as in \eqref{e:IL1}, we find
\begin{align} \nonumber
  (1+t)^{n/2}\|\mathcal{I}_1(t)\|_{L^\infty}
  &\le  C(1+t)^{n/2} \int_0^{t/2} \frac{1}{(t-s)^{n/2}}
  \|w(s)\|_{L^2}\|\Delta w(s)\|_{L^2}\,\rmd s\\ \label{e:IL2}
  &\le C\int_0^{t/2} \|w(s)\|_{L^2}\|\Delta w(s)\|_{L^2}\rmd s
  \le C \|w\|_Y^2\,.
\end{align}
On the other hand, using \eqref{e:lplq} with $(p,q) = (2,\infty)$
we obtain
\begin{align}\nonumber
  (1+t)^{n/2}\|\mathcal{I}_2(t)\|_{L^\infty}
  &\le  C(1+t)^{n/2} \int_{t/2}^t \frac{1}{(t-s)^{n/4}}
  \|w(s)\|_{L^\infty} \|\Delta w(s)\|_{L^2}\,\rmd s\\ \label{e:IL3}
  &\le C\|w\|_Y \frac{1}{(1+t)} \int_{t/2}^t
  \frac{1}{(t-s)^{n/4}}(1+s)\|\Delta w(s)\|_{L^2}\,\rmd s\\ \nonumber
  &\le C \|w\|_Y^2 \frac{1}{(1+t)^{\frac{n}4 +\frac1r}}\,,
\end{align}
where in the last line we used H\"older's inequality as in 
\eqref{e:Ismallt}. This shows that
\[
  \sup_{t\ge 0}(1+t)^{n/2}\|\mathcal{I}(t)\|_{L^\infty} \,\le\, 
  C \|w\|_Y^2\,.
\]

\noindent
\textbf{Estimate on $\|\Delta\mathcal{I}(t)\|_{L^2}$} \\
The estimates for the second derivative require maximal
regularity, Proposition~\ref{p:4}. We need to estimate 
$\int_0^\infty (1+t)^r\|\Delta\mathcal{I}(t)\|_{L^2}^r
\,\rmd t$. We therefore split the integral and first estimate
\[
  \int_0^T (1+t)^r\|\Delta\mathcal{I}(t)\|_{L^2}^r \,\rmd t 
  \,\le\, C\int_0^T \|w(t)\|_{L^\infty}^r 
  \|\Delta w(t)\|_{L^2}^r \rmd t \le C\|w\|_Y^{2r}\,,
\]
where we used Proposition~\ref{p:4} and estimate \eqref{e:nlest}
in the first inequality, and the definition of $\|w\|_Y$ in the 
second inequality.

We next derive estimates for $\Delta\mathcal{I}(t)$ at
$t\ge T$. Here it is more convenient to decompose
\begin{align*}
  \mathcal{I}(t) \,&=\, \int_0^{t-T}\mathcal{M}(t,s)
  K(s,w,\nabla w,\Delta w)\,\rmd s + 
  \int_{t-T}^t\mathcal{M}(t,s)K(s,w,\nabla w,\Delta w)
  \,\rmd s\\ 
  \,&=:\, \mathcal{I}_3(t) + \mathcal{I}_4(t)\,.
\end{align*}
Using \eqref{e:lplq} with $(p,q) = (1,2)$ and Remark~\ref{r:1}
we control the first term as
\begin{align*}
  (1+t)\|\Delta\mathcal{I}_3(t)\|_{L^2}
  &\,\le\, C(1+t) \int_0^{t-T}\frac{1}{(t-s)^{1+\frac{n}{4}}}
  \|w(s)\|_{L^2}\|\Delta w(s)\|_{L^2}\,\rmd s\\
  &\,\le\, C\|w\|_Y (1+t) \int_0^{t-T} \frac{1}{(t-s)^{1+\frac{n}{4}}} 
  \frac{1}{(1+s)^{1+\frac{n}{4}}}(1+s)\|\Delta w(s)\|_{L^2}\,\rmd s\\
  &\,\le\, \|w\|_Y^2 \frac{1}{(1+t)^{\frac{n}{4}}}\,,
\end{align*}
where in the last estimate we used H\"older's inequality together
with the fact that
\begin{equation}\label{e:auxbound}
  \left\{\int_0^t \Bigl(\frac{1}{1+t-s}\frac{1}{(t-s)^{n/4}} 
  \frac{1}{(1+s)^{1+\frac{n}{4}}}\Bigr)^q \,\rmd s\right\}^{1/q}
  \,\le\, \frac{C}{(1+t)^{1+\frac{n}{4}}}, \quad \hbox{for any }
  q \ge 1\,.
\end{equation}
Integrating over time and recalling that $nr/4 > 1$, we obtain 
the desired bound for $\mathcal{I}_3$:
\[
  \int_T^\infty (1+t)^r \|\Delta\mathcal{I}_3(t)\|_{L^2}^r
  \,\rmd t \,\le\, C\|w\|_Y^{2r}\,. 
\]
The other term $\mathcal{I}_4$ is estimated directly using 
\eqref{e:nlest} and Corollary~\ref{c:max}:
\[
  \int_T^\infty (1+t)^r \|\Delta\mathcal{I}_4(t)\|_{L^2}^r\,
  \rmd t \,\le\, C\int_0^\infty (1+t)^r \|w(t)\|_{L^\infty}^r
  \|\Delta w(t)\|_{L^2}^r \,\rmd t \,\le\, C\|w\|_Y^{2r}\,.
\]

\noindent
\textbf{Estimates on $\|\mathcal{J}(t)\|_{L^1}$ and
$\|\mathcal{J}(t)\|_{L^\infty}$} \\
In this term, the nonlinearity does not contain derivatives 
which would yield decay, but we can exploit the the stronger 
decay of the linear operator $\mathcal{M}(t,s)$ when acting
on $H(s,w)$. Indeed, since $H(s,w) = (0,\hat g_2(s,w))^T$, 
we have
\[
  \mathcal{M}(t,s) H(s,w) \,=\, \left(\begin{array}{l}
  \mathcal{M}_{0\mathrm{h}}(t,s)\hat g_2(s,w) \\
  \mathcal{M}_{\mathrm{hh}}(t,s)\hat g_2(s,w)\end{array}\right)
  \,=:\, \mathcal{M}_{\cdot\mathrm{h}}(t,s)\hat g_2(s,w)\,,
\]
and it follows from Proposition~\ref{p:3} that
\begin{equation}\label{e:MHest}
  \|\mathcal{M}(t,s)H(s,w)\|_{L^q(\R^n)} \,\le\, C\,\frac{1}{1+t-s}
  \frac{1}{(t-s)^{\frac{n}2(\frac1p-\frac1q)}}\ \|H(s,w)\|_{L^p(\R^n)}\,,
\end{equation}
for $1 \le p \le q \le \infty$. Using \eqref{e:MHest} with $p = q = 1$ 
and estimate \eqref{e:nlesth}, we thus find
\begin{equation}\label{e:Jest1}
  \|\mathcal{J}(t)\|_{L^1} \,\le\, C\int_0^t \frac{1}{1+t-s}
   \|w(s)\|_{L^2}^2\,\rmd s \,\le\, C\|w\|_Y^2 \int_0^t \frac{1}{1+t-s} 
  \frac{1}{(1+s)^{n/2}}\,\rmd s \,\le\, C \|w\|_Y^2\,.
\end{equation}
In a similar way, using \eqref{e:MHest} with $(p,q) = (2,\infty)$, 
we arrive at
\begin{align}\nonumber
  (1+t)^{n/2}\|\mathcal{J}(t)\|_{L^\infty} 
  \,&\le\, C(1+t)^{n/2} \int_0^t \frac{1}{1+t-s}\frac{1}{(t-s)^{n/4}} 
  \|w(s)\|_{L^\infty} \|w(s)\|_{L^2}\,\rmd s\\ \label{e:Jest2}
  \,&\le\, C \|w\|_Y^2 (1+t)^{n/2} \int_0^t \frac{1}{1+t-s}\frac{1}
  {(t-s)^{n/4}} \frac{1}{(1+s)^{3n/4}}\,\rmd s\\ \nonumber
  \,&\le\, C\|w\|_Y^2\,.
\end{align}

\noindent
\textbf{Estimate on $\|\Delta\mathcal{J}(t)\|_{L^2}$} \\
We first observe that 
\[
  \Delta\mathcal{J}(t) \,=\, \int_0^t \mathcal{M}_{\cdot\mathrm{h}}
  (t,s)\Delta \hat g_2(s,w)\,\rmd s\,.
\]
Since $\hat g_2(t,w) = g_2(v_*(t)+w)[w_\mathrm{h},w_\mathrm{h}]$, 
it is clear that $|\Delta \hat g_2(t,w)| \le C (|w||\Delta w|
+|\nabla w|^2)$, so that $\Delta \hat g_2$ satisfies the
same estimates \eqref{e:nlest} as $F$ and $G$. Using the 
$L^1$--$L^2$ estimate for $\mathcal{M}_{\cdot\mathrm{h}}(t,s)$, 
we thus find
\begin{align*}
  (1+t)\|\Delta \mathcal{J}(t)\|_{L^2} 
  \,&\le\, C(1+t) \int_0^t \frac{1}{1+t-s}\frac{1}{(t-s)^{n/4}} 
  \|w(s)\|_{L^2} \|\Delta w(s)\|_{L^2}\,\rmd s\\
  \,&\le\, C\|w\|_Y (1+t) \int_0^t \frac{1}{1+t-s}\frac{1}{(t-s)^{n/4}} 
  \frac{1}{(1+s)^{1+\frac{n}{4}}}(1+s)\|\Delta w(s)\|_{L^2}\,\rmd s\\
  \,&\le\, C\|w\|_Y^2 \frac{1}{(1+t)^{n/4}}\,,
\end{align*}
where in the last line we used H\"older's inequality and estimate 
\eqref{e:auxbound}. We deduce that
\[
  \int_0^\infty (1+t)^r \|\Delta \mathcal{J}(t)\|_{L^2}^r
  \,\rmd t \,\le\, C\|w\|_Y^{2r}\,, 
\]
and the proof of Proposition~\ref{p:fp} is complete. 
\end{Proof}

\noindent
\textbf{Proof of Theorem~\ref{t:1}} ($n \le 3$). As was observed in
Section~\ref{s:nf}, we can work with the transformed
equation~\eqref{e:ednfp} instead of the original perturbation equation
\eqref{e:pertv}. Also, we can assume without loss of generality that
the initial perturbation $w_0$ satisfies $\|w_0\|_{L^1} +
\|w_0\|_{H^2} \le \delta_0$ for some small $\delta_0 > 0$.  Under
these assumptions, we can solve equation~\eqref{e:mild} by a standard
fixed point argument in the Banach space $Y$ defined by
\eqref{e:norm}. Indeed, let $\mathcal{N}$ denote the right-hand 
side of \eqref{e:mild}, namely
\[
  (\mathcal{N}w)(t) \,=\,\mathcal{M}(t,0)w_0 + \int_0^t \mathcal{M}(t,s)
  K(s,w,\nabla w,\Delta w) \,\rmd s + \int_0^t \mathcal{M}(t,s) 
  H(s,w) \,\rmd s\,. 
\]
If $w \in Y$ satisfies $\|w\|_Y \le \delta_2$, where $\delta_2 > 0$
is as in Proposition~\ref{p:fp}, we know that $\mathcal{N}w \in Y$
and that
\begin{equation}\label{e:NN1}
  \|\mathcal{N}w\|_Y \,\le\, C_1 \delta_0 + C_2 \|w\|_Y^2\,.
\end{equation}
Similar calculations show that
\begin{equation}\label{e:NN2}
  \|\mathcal{N}w - \mathcal{N}\tilde w\|_Y \,\le\, C_2 
  (\|w\|_Y + \|\tilde w\|_Y)\|w - \tilde w\|_Y\,,
\end{equation}
whenever $w, \tilde w \in Y$ with $\|w\|_Y \le \delta_2$, $\|\tilde
w\|_Y \le \delta_2$. Let $\mathcal{B} \subset Y$ denotes the ball of
radius $R = \min(2C_1\delta_0,\delta_2)$ centered at the origin. If
$\delta_0 > 0$ is small enough so that $2C_2 R < 1$, it follows easily
from \eqref{e:NN1}, \eqref{e:NN2} that $\mathcal{N}(\mathcal{B})
\subset \mathcal{B}$ and that $\mathcal{N}$ is a strict contraction in
$\mathcal{B}$. Let $w \in Y$ be the unique fixed point of
$\mathcal{N}$ in $\mathcal{B}$. Then $w$ is a global solution of 
\eqref{e:ednfp}, and if we return to the original variables by 
setting $u(t,x) = \Psi(v_*(t) + w(t,x))$, we obtain a global solution
of \eqref{rdsystem} which satisfies the decay estimate
\eqref{e:t1decay} (with $t_0 = 0$), 
because
\[
  |u(t,x) - u_*(t)| \,=\, |\Psi(v_*(t) + w(t,x)) - \Psi(v_*(t))|
  \,\le\, C\,\frac{\|w\|_Y}{(1+t)^{n/2}}\,,\quad t \ge 0\,.
\]
This concludes the proof. \hfill $\Box$

\begin{Remark}\label{r:3}
The limitation $n \le 3$ in the above proof is due to the fact that
we use maximal regularity (MR) in the Hilbert space $L^2(\R^n)$ only, see
e.g. \eqref{e:Ismallt}. This choice was made for simplicity, but for
equation \eqref{e:ednfls} it is known that MR holds in all $L^p$
spaces with $1 < p < \infty$, see \cite{HM00,PS01}. It is not 
difficult to verify that the argument above can be adapted to any 
space dimension $n$ if we use MR in $L^p(\R^n)$ with $p$ sufficiently 
large, depending on $n$.
\end{Remark}

\section{Asymptotic behavior}\label{s:asym}

We know from Theorem~\ref{t:1} that small, localized perturbations
of the periodic solution $u_*(t)$ converge to zero like $t^{-n/2}$
as $t \to +\infty$. This decay rate is optimal in general, and it 
is even possible to compute the leading term in the asymptotic 
expansion of the perturbation as $t \to +\infty$. In this section, 
we assume (for simplicity) that $1 \le n \le 3$ and we consider 
the solution $w(t,x)$ of \eqref{e:ednfp} with small initial 
data $w_0 \in L^1(\R^n) \cap H^2(\R^n)$. If we decompose this 
solution as $w(t,x) = (w_0(t,x),w_{\mathrm{h}}(t,x))^T$, we first
observe that the hyperbolic part $w_{\mathrm{h}}(t,x) \in \R^{N-1}$ 
has a faster decay as $t \to \infty$. 
 
\begin{Proposition}\label{p:faster}
If the initial data $w_0 \in L^1(\R^n) \cap H^2(\R^n)$ are
sufficiently small, the hyperbolic component of the solution 
$w$ of \eqref{e:ednfp} satisfies
\[
  \sup_{t \ge 0} (1+t)\|w_{\mathrm{h}}(t)\|_{L^1} + \sup_{t \ge 0}
  (1+t)^{1+\frac{n}2} \|w_{\mathrm{h}}(t)\|_{L^\infty} \,\le\, 
  C(\|w_0\|_{L^1} + \|w_0\|_{H^2})\,. 
\]
\end{Proposition}

\begin{Proof}
Projecting the integral equation \eqref{e:mild} onto the 
hyperbolic component, we find
\begin{align*}
  w_\mathrm{h}(t)\,&=\,\mathcal{M}_{\mathrm{h}\cdot}(t,0)w_0
  + \int_0^t \mathcal{M}_{\mathrm{h}\cdot}(t,s)K(s,w,
  \nabla w,\Delta w)\,\rmd s + 
  \int_0^t \mathcal{M}_{\mathrm{h\cdot}}(t,s)H(s,w)\,\rmd s \\
  \,&=:\,\mathcal{M}_{\mathrm{h}\cdot}(t,0)w_0 + 
  \mathcal{I}_{\mathrm{h}}(t) + \mathcal{J}_{\mathrm{h}}(t)\,,
\end{align*}
where $\mathcal{M}_{\mathrm{h\cdot}}(t,s) = (\mathcal{M}_{
\mathrm{h0}}(t,s),\mathcal{M}_{\mathrm{hh}}(t,s))$. We know from
Proposition~\ref{p:3} that 
\begin{equation}\label{e:MHest2}
  \|\mathcal{M}_{\mathrm{h\cdot}}(t,s)w\|_{L^q(\R^n)} \le C\,\frac{1}{1+t-s}
  \frac{1}{(t-s)^{\frac{n}2(\frac1p-\frac1q)}}\ \|w\|_{L^p(\R^n)}\,,
\end{equation}
for $1 \le p \le q \le \infty$. In particular, we have
\[
  \sup_{t \ge 0} (1+t)\|\mathcal{M}_{\mathrm{h}\cdot}(t,0)w_0\|_{L^1} 
  + \sup_{t \ge 0} (1+t)^{1+\frac{n}2} \|\mathcal{M}_{\mathrm{h}\cdot}(t,0)
  w_0 \|_{L^\infty} \,\le\, C(\|w_0\|_{L^1} + \|w_0\|_{L^\infty})\,. 
\]
Moreover, proceeding as in the proof of Proposition~\ref{p:fp}
we obtain
\[
  \|\mathcal{I}_{\mathrm{h}}(t)\|_{L^1} \,\le\, C\int_0^t 
  \frac{1}{1+t-s}\frac{\|w\|_Y}{(1+s)^{1+\frac{n}{4}}}
  \,(1+s)\|\Delta w(s)\|_{L^2}\,\rmd s \,\le\, \frac{C}{1+t}
  \,\|w\|_Y^2\,,
\]
and the same result holds for $(1+t)^{n/2}\|\mathcal{I}_{
\mathrm{h}}(t)\|_{L^\infty}$. Finally, to bound the term 
$\mathcal{J}_{\mathrm{h}}$, we observe that 
$\mathcal{M}_{\mathrm{h\cdot}}(t,s)H(s,w) = \mathcal{M}_{
\mathrm{hh}}(t,s)\hat g_2(s,w)$ and we use the strong decay in 
time given by Proposition~\ref{p:3}. We thus find
\begin{align}\label{e:Jhest}
  \|\mathcal{J}_{\mathrm{h}}(t)\|_{L^1} \,&\le\, C\int_0^t 
  \frac{1}{(1+t-s)^2}\,\|w_{\mathrm{h}}(s)\|_{L^1}
  \|w_{\mathrm{h}}(s)\|_{L^\infty}\,\rmd s \\ \nonumber
  \,&\le\, C\int_0^t \frac{1}{(1+t-s)^2}\frac{\|w\|_Y^2}{(1+s)^{n/2}}
  \,\rmd s \,\le\, \frac{C}{(1+t)^{n/2}}\,\|w\|_Y^2\,,
\end{align}
and the same result holds for $(1+t)^{n/2}\|\mathcal{J}_{
\mathrm{h}}(t)\|_{L^\infty}$. This gives the desired result if $n 
\ge 2$. If $n = 1$, we only have
\[
  \sup_{t \ge 0} (1+t)^{1/2}\|w_{\mathrm{h}}(t)\|_{L^1} + \sup_{t \ge 0}
  (1+t) \|w_{\mathrm{h}}(t)\|_{L^\infty} \,\le\, C\|w_0\|_{L^1\cap H^2}\,,
\]
but if we now return to \eqref{e:Jhest} we obtain the stronger
estimate
\[
  \|\mathcal{J}_{\mathrm{h}}(t)\|_{L^1} \,\le\, C\|w_0\|_{L^1\cap H^2}^2
  \int_0^t \frac{1}{(1+t-s)^2}\,\frac{1}{(1+s)^{3/2}}\,\rmd s
  \,\le\, C\,\frac{\|w_0\|_{L^1\cap H^2}^2}{(1+t)^{3/2}}\,, 
\]
which also holds for $(1+t)^{1/2}\|\mathcal{J}_{\mathrm{h}}(t)
\|_{L^\infty}$. This concludes the proof. 
\end{Proof}

We next consider the central component $w_0(t,x) \in \R$, and prove 
that it behaves asymptotically like a solution of a linear equation 
with suitably modified initial data. 

\begin{Proposition}\label{p:central}
If the initial data $w_0 \in L^1(\R^n) \cap H^2(\R^n)$ are
sufficiently small, the central component of the solution 
$w$ of \eqref{e:ednfp} satisfies
\[
  \|w_0(t) - \mathcal{M}_{00}(t,0)w_\infty\|_{L^1} + (1+t)^{n/2} 
  \|w_0(t) - \mathcal{M}_{00}(t,0)w_\infty\|_{L^\infty} 
  \,\le\, \frac{C}{(1+t)^\gamma}\,(\|w_0\|_{L^1} + \|w_0\|_{H^2})\,,
\]
where $\gamma = \frac{n}{4} + \frac{1}{r} < 1$ and $w_\infty \in 
L^1(\R^n) \cap L^\infty(\R^n)$ is defined by
\[
  w_\infty \,=\, (w_0)_0 + \int_0^\infty K_0(s,w(s),\nabla w(s),
  \Delta w(s))\,\rmd s\,.
\]
\end{Proposition}

\begin{Proof}
Projecting \eqref{e:mild} onto the central component, we find
\begin{align*}
  w_0(t)\,&=\,\mathcal{M}_{0\cdot}(t,0)w_0 + \int_0^t 
  \mathcal{M}_{0\cdot}(t,s)K(s,w,\nabla w,\Delta w)\,\rmd s + 
  \int_0^t \mathcal{M}_{0\cdot}(t,s)H(s,w)\,\rmd s \\
  \,&=:\,\mathcal{M}_{0\cdot}(t,0)w_0 + \mathcal{I}_0(t) + 
  \mathcal{J}_0(t)\,,
\end{align*}
where $\mathcal{M}_{0\cdot}(t,s) = (\mathcal{M}_{00}(t,s),
\mathcal{M}_{0\mathrm{h}}(t,s))$. Our goal is to extract from 
$w_0(t)$ the leading contributions as $t \to +\infty$. The 
last term $\mathcal{J}_0(t)$ is clearly negligible in this 
limit. Indeed, using Proposition~\ref{p:faster} and proceeding 
as in \eqref{e:Jest1}, \eqref{e:Jest2}, we obtain
\begin{equation}\label{e:negl}
  \|\mathcal{J}_0(t)\|_{L^1} + (1+t)^{n/2} \|\mathcal{J}_0(t)
  \|_{L^\infty} \,\le\, \frac{C}{1+t}\,\|w_0\|_{L^1\cap H^2}\,.
\end{equation}
The same estimate holds for the linear term $\mathcal{M}_{0\mathrm{h}}
(t,0)(w_0)_\mathrm{h}$, because $\mathcal{M}_{0\mathrm{h}}(t,0)$
decays as fast as $(1+t)^{-1}\mathcal{M}_{00}(t,0)$ by 
Proposition~\ref{p:3}. Using the same remark and proceeding 
as in \eqref{e:IL1}, \eqref{e:IL2}, \eqref{e:IL3}, we see that
the integral term $\mathcal{I}_{0\mathrm{h}}(t) := \int_0^t 
\mathcal{M}_{0\mathrm{h}}(t,s)K_\mathrm{h}(s,w,\nabla w,\Delta w)
\,\rmd s$ also satisfies \eqref{e:negl}. So the only remaining 
terms are $\mathcal{M}_{00}(t,0)(w_0)_0$ and
\begin{align*}
  \mathcal{I}_{00}(t) \,&=\, \int_0^t \mathcal{M}_{00}(t,s)
  K_0(s,w,\nabla w,\Delta w)\,\rmd s \\
  \,&=\, \int_{t/2}^t \mathcal{M}_{00}(t,s) K_0(s,w,\nabla w,
    \Delta w)\,\rmd s 
  \,+\, \int_0^{t/2} (\mathcal{M}_{00}(t,s)-\mathcal{M}_{00}(t,0))
     K_0(s,w,\nabla w,\Delta w)\,\rmd s \\ 
  &\quad + \mathcal{M}_{00}(t,0)\int_0^\infty K_0(s,w,\nabla w,
    \Delta w)\,\rmd s 
  \,-\, \mathcal{M}_{00}(t,0)\int_{t/2}^\infty K_0(s,w,\nabla w,
    \Delta w)\,\rmd s \\ 
  \,&=:\, \mathcal{I}_{01}(t) + \mathcal{I}_{02}(t) + 
    \mathcal{I}_{03}(t) + \mathcal{I}_{04}(t)\,.
\end{align*}
Proceeding as in \eqref{e:IL1}, \eqref{e:IL3}, it is straightforward 
to verify that 
\[
  \|\mathcal{I}_{01}(t)\|_{L^1} + (1+t)^{n/2} \|\mathcal{I}_{01}(t)
  \|_{L^\infty} \,\le\, \frac{C}{(1+t)^\gamma}\,\|w\|_Y^2\,,
  \quad t \ge 0\,,
\]
where $\gamma = \frac{n}{4} + \frac{1}{r} < 1$, and the same estimate 
clearly holds for $\mathcal{I}_{04}(t)$ too. Finally, using 
Remark~\ref{r:2} to bound the difference $\mathcal{M}_{00}(t,s)-
\mathcal{M}_{00}(t,0)$, we obtain
\[
  \|\mathcal{I}_{02}(t)\|_{L^1} + (1+t)^{n/2} \|\mathcal{I}_{02}(t)
  \|_{L^\infty} \,\le\, C \int_0^{t/2} \frac{s}{t}\,
  \|K_0(s,w,\nabla w,\Delta w)\|_{L^1}\,\rmd s \,\le\, 
  \frac{C}{(1+t)^\gamma}\,\|w\|_Y^2\,.
\]
This concludes the proof of Proposition~\ref{p:central}, because 
$\mathcal{M}_{00}(t,0)(w_0)_0 + \mathcal{I}_{03}(t) = 
\mathcal{M}_{00}(t,0)w_\infty$.
\end{Proof}

It is now rather easy to prove Theorem~\ref{t:2} (in the case where 
$n \le 3$). Combining Propositions~\ref{p:faster}
and \ref{p:central}, we find
\begin{equation}\label{e:fin1}
  \|w(t) - e_1 W(t)\|_{L^1} + (1+t)^{n/2}\|w(t) - e_1
  W(t)\|_{L^\infty} \,\le\, \frac{C}{(1+t)^\gamma}\,(\|w_0\|_{L^1} 
  + \|w_0\|_{H^2})\,,
\end{equation}
where $W(t) = \mathcal{M}_{00}(t,0)w_\infty$ and $e_1 = (1,0)^T$ 
is the first vector of the canonical basis in $\R^N$. Furthermore, 
we claim that 
\begin{equation}\label{e:fin2}
  \|t^{n/2}W(t,xt^{1/2}) - \tilde \alpha G\|_{L^1 \cap 
  L^\infty} ~\xrightarrow[t\to +\infty]{}~ 0\,,
\end{equation}
where $G$ is defined in \eqref{e:GGdef} and
\begin{equation}\label{e:fin3}
  \tilde \alpha \,=\, \int_{\R^n} w_\infty(x)\,\rmd x \,=\,
  e_1^T \left(\int_{\R^n} w_0(x)\,\rmd x + \int_0^\infty \!\!\int_{\R^n}
  K(t,w,\nabla w,\Delta w)\,\rmd x\,\rmd t\right)\,.   
\end{equation}
To prove the $L^\infty$ claim in \eqref{e:fin2}, we use Fourier 
transforms and simply note that the quantity 
\[
  \|\hat W(t,kt^{-1/2}) - \tilde \alpha \hat{G}(k)\|_{L^1}
  \,=\, \|M_{00}(t,0;kt^{-1/2})\hat w_\infty(k t^{-1/2}) - e^{-d_0 k^2}
  \hat w_\infty(0)\|_{L^1}
\]
converges to zero as $t \to \infty$ by Lebesgue's dominated
convergence theorem, in view of Proposition~\ref{p:1} and 
Remark~\ref{r:asym}. The $L^1$ claim can be established in 
a similar way, using the same ideas as in the proof of 
Proposition~\ref{p:3} (we omit the details). 

We now return to the original variables. Since $u_*(t) = 
\Psi(\omega t e_1)$, the solution of \eqref{rdsystem} 
given by $u(t,x) = \Psi(v(t,x)) = \Psi(v_*(t) + w(t,x))$ 
can be decomposed as in \eqref{e:udecomp}, with $\alpha(t,x) = 
\omega^{-1}W(t,x)$ and
\[
  \beta(t,x) \,=\, \Psi(v_*(t) + w(t,x)) - \Psi(v_*(t)) - 
  \Psi'(v_*(t))w(t,x) + \Psi'(v_*(t))(w(t,x)-e_1 W(t,x))\,.
\]
Estimates \eqref{e:fin1}, \eqref{e:fin2} immediately give
\eqref{e:betaest}, \eqref{e:alphaest} with $\alpha_* =
\omega^{-1}\tilde \alpha$. Finally, the formula \eqref{e:alphastardef}
for $\alpha_*$ follows from the expression \eqref{e:fin3} of $\tilde
\alpha$ and the fact that $U_*(0) = (\Psi'(0)^{-1})^T e_1$, $u_*'(0) =
\omega \Psi'(0)e_1$. This concludes the proof of Theorem~\ref{t:2}.
\hfill $\Box$

\section{Examples and perspectives}\label{s:ex}

In this final section, we first give a simple example of 
a 2-species reaction-diffusion system with a periodic orbit
$u_*(t)$ which is asymptotically stable for the ODE dynamics 
but does not satisfy Hypothesis~\ref{h:2}. We then discuss 
possible generalizations of the results of this paper. 

\subsection{Destabilization by diffusion: a simple example}

One may feel inclined to believe that ODE-stable periodic orbits tend
to be stable for the PDE dynamics, that is, that our Hypothesis
\ref{h:2} is satisfied in most cases where the periodic orbit is
stable for the ODE. Our example below shows that this is not the case,
even for a simple reaction-diffusion system with only two species.

We consider the following $2$-species reaction-diffusion system
\begin{equation}\label{e:ex1}
  u_t \,=\, D u_{xx} + Ju + (\epsilon^2-|u|^2)Ru \,,
\end{equation}
where $u = (u_1,u_2)^T \in \R^2$ and $|u|^2 = u_1^2 + u_2^2$. 
Here $\epsilon > 0$ is a parameter, $D$ is a $2\times2$ real matrix 
with positive eigenvalues, and
\[
  J \,=\, \begin{pmatrix} 0 & -1 \\ 1 & 0 \end{pmatrix}\,, \qquad
  R \,=\, \begin{pmatrix}  \cos(\theta) & -\sin(\theta) \\
  \sin(\theta) & \cos(\theta)\end{pmatrix}\,, \qquad \theta \in 
  (-\pi/2,\pi/2)\,.
\]
The system~\eqref{e:ex1} has a $2\pi$-periodic solution $u_*(t) = \epsilon 
\bar u(t)$, where $\bar u(t) = (\cos(t),\sin(t))^T$. Linearizing 
\eqref{e:ex1} at $u_*(t)$ we obtain
\[
  v_t \,=\, D v_{xx} + (J - 2\epsilon^2 R\,\bar u(t)\,\bar u(t)^T)v\,,
\]
or equivalently
\begin{equation}\label{e:ex2}
  v_t \,=\, -k^2 D v + (J - 2\epsilon^2 R\,\bar u(t)\,\bar u(t)^T)v\,,
  \quad k \in \R\,.
\end{equation}
Of course $v(t) = \bar u'(t)$ is a solution of \eqref{e:ex2} for 
$k = 0$.
 
Let $F = F_0(2\pi,0)$ be the Floquet matrix associated to~\eqref{e:ex2} for $k = 0$. Then
\[
  \Det(F) \,=\, \exp\Bigl(\Tr \int_0^{2\pi} (J - 2\epsilon^2 
  R\,\bar u(t)\,\bar u(t)^T)\,\rmd t\Bigr) \,=\, \exp(-2\pi
  \epsilon^2\Tr(R))\,.
\]
The Floquet exponents (for $k = 0$) are therefore $\lambda_1 = 0$
and $\lambda_2 = -\epsilon^2 \Tr(R)$. As $\Tr(R) = 2\cos(\theta) > 0$, 
it follows that $u_*(t)$ is a stable periodic orbit for the ODE 
dynamics associated to \eqref{e:ex1}.

To compute the Floquet exponents for small $k$, we consider the 
adjoint ODE
\begin{equation}\label{e:ex3}
  U_t \,=\, (J + 2 \epsilon^2 \bar u(t)\,\bar u(t)^T R^T)U\,.
\end{equation}
As is easily verified, the unique nontrivial bounded solution of 
\eqref{e:ex3} is $U_*(t) = R \bar u'(t)$. Using formula~\eqref{e:d0def}, 
we conclude that  
\[
  d_0 \,=\, \frac{\int_0^{2\pi} \bar u'(t)^T R^T D \,\bar u'(t)
  \,\rmd t}{\int_0^{2\pi} \bar u'(t)^T R^T \bar u'(t)\,\rmd t} 
  \,=\, \frac{\Tr(R^T D)}{\Tr(R^T)} \,=\, \frac12 \Bigl(
  \Tr(D) - \tan(\theta)\Tr(JD)\Bigr)\,.
\]
If the diffusion matrix $D$ is symmetric, then $\Tr(D) > 0$ and
$\Tr(JD) = 0$, hence necessarily $d_0 > 0$, which means that the
periodic solution $u_*(t)$ is spectrally stable for long wave-length
perturbations. But if $D$ is a nonsymmetric matrix, then $\Tr(JD) \neq
0$ and therefore we can choose $\theta \in (-\pi/2,\pi/2)$ in such a
way that $d_0 < 0$. This gives an example of a periodic orbit
exhibiting a {\em sideband} instability.

On the other hand, as $\epsilon \to 0$, the periodic orbit $u_*(t)$
reduces to the fixed point $u = 0$, for which it is easy to perform a
stability analysis. For a fixed wavenumber $k \in \R$, we have to
compute the eigenvalues $\lambda_1(k), \lambda_2(k)$ of the linearized
operator $J-Dk^2$. By direct calculation, we find
\[
  \Tr(J-Dk^2) \,=\, -\Tr(D)k^2 \,\le\, 0\,, \quad \hbox{and} \quad
  \Det(J-Dk^2) \,=\, 1 + \Tr(JD)k^2 + \Det(D)k^4\,.
\]
If we choose $D$ such that $\Tr(JD) + 2(\Det(D))^{1/2} < 0$, we see
that there exists a nonempty open interval $I \subset (0,+\infty)$ such
that $\Tr(J-Dk^2) < 0$ and $\Det(J-Dk^2) < 0$ if $k^2 \in I$.  Thus,
one of the eigenvalues $\lambda_i(k)$ is strictly positive, which means
that the equilibrium $u = 0$ is unstable with respect to perturbations
with wavenumbers $k$ such that $k^2 \in I$.  By continuity, this
Turing instability persists for the periodic orbit $u_*(t) = \epsilon
\bar u(t)$ if $\epsilon > 0$ is sufficiently small: for $k^2 \in I$,
one of the Floquet exponents has positive real part.

Summarizing, the periodic solution $u_*(t) = \epsilon \bar u(t)$ 
of the reaction-diffusion system~\eqref{e:ex1} exhibits:
\vspace{-2mm}
\begin{enumerate}
\setlength{\itemsep}{-0.4mm}
\item a {\em sideband instability}, if $\Tr(D) - \tan(\theta)
\Tr(JD) < 0$;
\item a {\em Turing instability}, if $\Tr(JD) + 2(\Det(D))^{1/2} 
< 0$ and $\epsilon \ll 1$. 
\end{enumerate}

\begin{Remarks} \hspace{1mm}\\[1mm]
{\bf 1.} 
The instability criteria above are never satisfied if the 
matrix $D$ is symmetric. But if $D$ has eigenvalues $d_1 > d_2 > 0$,
we can choose an invertible matrix $S$ so that $S^{-1}DS = 
\mathcal{D} = \mathrm{diag}(d_1,d_2)$. Then setting $u = Sw$ we 
obtain the equivalent system
\begin{equation}\label{e:ex4}
  w_t \,=\, \mathcal{D} w_{xx} + S^{-1}JSw + S^{-1}(1-|Sw|^2)RSw\,,
\end{equation}
where now the diffusion matrix has the usual, diagonal form. 
The Floquet exponents characterizing the stability properties of 
the periodic orbit are of course unaffected by this linear
transformation. Thus we can find 2-species systems of the 
form \eqref{e:ex4}, with diagonal diffusion matrix, which exhibit
either a sideband or a Turing instability.

\noindent{\bf 2.}
Our example is clearly reminiscent of the complex Ginzburg-Landau 
equation (CGLE),
\begin{equation}\label{e:CGL}
  u_t \,=\, (1+ia)\Delta u + u - (1-ic)|u|^2 u\,,
\end{equation}
where $a,c$ are real parameters and $u : \R_+ \times \R^n \to \C$.
The system~\eqref{e:CGL} possesses a homogeneous time-periodic solution of the form
$u(t,x) = \rme^{ict}$, which exhibits a sideband instability if $ac >
1$ ({\em Benjamin-Feir criterion}) and a Turing instability in other
parameter regions, see e.g.  \cite{AK02}. The complex Ginzburg-Landau
equation arises as a modulation equation near Hopf bifurcations in
reaction-diffusion systems (see, for example, \cite{sh,Mi}). One
therefore expects stability and instability properties of
small-amplitude periodic solutions near Hopf bifurcation to be
governed by those of the CGLE; see, for example, \cite{Ri2} for a
result in this direction.  Since the diffusion matrix of CGLE
possesses complex eigenvalues, it cannot be cast as a real
reaction-diffusion system with diagonal diffusion matrix. Our example
and the remark above show that one can almost explicitly recover the
properties of CGLE with diagonal diffusion matrices. Upon substituting
the diffusion matrix of CGLE in our example, one would recover
precisely the Benjamin-Feir criterion from the instability criterion
$\Tr(D) - \tan(\theta)\Tr(JD) < 0$. In a different direction, one
could also extend our results to diffusion matrices with $D+D^T>0$
without any additional difficulties, which would then include CGLE as
a particular example. In the context of reaction-diffusion modeling,
cross-diffusion phenomena that are associated with off-diagonal
elements of $D$ are however quite uncommon.
\end{Remarks}

\subsection{Discussion and perspectives}

We believe that the method presented here can be adapted to other
situations. We mention the stability of wave trains, $u(kx-\omega t)$,
with $u(\xi)=u(\xi+2\pi)$ and $\omega,k>0$, and Turing patterns
$u(kx)=u(-kx)$, with $u(\xi)=u(\xi+2\pi)$ and $k>0$. In both cases, one
finds continuous spectrum with diffusive decay properties for the
linearization. In both cases, the absence of relevant, self-coupling
terms in the neutral mode has been shown previously; see
\cite{Sc96,DSSS}.

Interesting questions arise when one attempts to extend the class of
allowed perturbations. One may for instance consider perturbations $v$
such that $\nabla v\in L^1$. In one space-dimension, this would
correspond to perturbations with different phase shifts at
$x=\pm\infty$. One would still expect the diffusive linear part to be
dominant so that one would find error function asymptotics for the
phase correction.

In fact, one would expect some type of stability for much more general
perturbations. For instance, in one space-dimension, the homogeneous
oscillation $u_*(t)$ is embedded in a family of wave train solutions
$u(kx-\omega t;k)$, with $k\approx 0$ and $u\approx u_*$;
\cite[Section~3.3]{ssdefect}.  Under our stability assumptions, using
Lyapunov-Schmidt reduction, one finds $\omega=\omega(k) =  \omega_0 + \omega_2
k^2+\rmO(k^4)$; see for instance \cite[Lemma 2.1]{RS}. Of course,
$u(kx-\omega t;k)$ is not close to any fixed homogeneous oscillation,
but it is close to an appropriate phase shift of the oscillation in
any finite region of space. All solutions in the basin of attraction of
such wave trains then stay close to our homogeneous oscillation,
orbitally, and pointwise in space. A natural question then asks for
the asymptotics of initial conditions of the type $u(k(x)
x-\omega(k(x)) t;k(x))$, where $k(x)\to k_\pm$ for $x\to\pm\infty$. In
the case of the real Ginzburg-Landau equation, which exhibits spatial
oscillations $u(kx)$, this question was addressed in \cite{CEE92,GM},
showing that asymptotics are governed by a nonlinear diffusion
equation $\theta_t=W(\theta_x)_x$, with locally uniform convergence to
a fixed, \emph{intermediate} wavenumber.  Near temporal oscillations, one
expects dynamics to be governed by a viscous conservation law
$\theta_t=d(\theta_x)_x+j(\theta_x)$, so that solutions with
asymptotically constant wavenumber would be expected to converge to
viscous shocks \cite[Section 8]{DSSS}, or to rarefaction waves.

More generally, one could ask about the orbital stability of a family
of oscillations: starting with an initial condition $u(k(x)
x-\omega(k(x)) t;k(x))$, with $k(x)$ bounded and $k'$ small, will the
solutions remain close to a solution of that form for all times?
Again, this question possesses a simple answer for the approximation
of the dynamics by a viscous conservation law due to the maximum
principle, which gives immediate supremum bounds on $\theta_x$ in
terms of the initial condition. For the full reaction-diffusion
dynamics however, which are only approximately described by a viscous
conservation law, this question remains wide open.

\end{document}